\documentclass[12pt,reqno]{amsart}
\usepackage{amssymb,amsmath,tabularx,mathrsfs,mathbbol,yfonts,upgreek,hyperref}
\usepackage{amsthm,verbatim,comment}
\usepackage{geometry}
\usepackage{stmaryrd}
\usepackage{paralist}
\usepackage[all]{xy}
\usepackage{mathdots}
\usepackage{tikz}
\usepackage{ytableau}
\usetikzlibrary{arrows,matrix,positioning,fit,patterns,calc,decorations.pathreplacing}
\usepackage{graphicx}
\allowdisplaybreaks
\usepackage{nicematrix}
\usepackage[capitalise]{cleveref}
\usepackage{youngtab}

\pgfdeclarelayer{bg}    

\newdimen\leftpt
\newdimen\righpt
\tikzset{%
  pics/overbrace/.style args={#1,#2,#3}{%
    code = {
       \pgfextractx\leftpt{\pgfpointanchor{m-1-#1}{north west}}  
       \pgfextractx\righpt{\pgfpointanchor{m-1-#2}{north east}}
       \node[rectangle, above delimiter={\{}, minimum width=\the\dimexpr\righpt-\leftpt+1pt,
             label={[yshift=3mm]above:$ #3$}]
             at ($ (m-1-#1)!0.5!(m-1-#2)+(0,0.2) $){};
       }
   },
   pics/leftbrace/.style args={#1,#2,#3}{%
   code = {
       \pgfextracty\leftpt{\pgfpointanchor{m-#1-1}{north west}}  
       \pgfextracty\righpt{\pgfpointanchor{m-#2-1}{south west}}
       \node[rectangle, left delimiter={\{}, minimum height=\dimexpr\leftpt-\righpt,
             label={[xshift=-3mm]left:$#3 $}]
             at ($ (m-#1-1)!0.45!(m-#2-1)+(0,0) $){};
       }
   }
}

\NiceMatrixOptions{
code-for-first-row = \color{red} ,
code-for-last-row = \color{red} ,
code-for-first-col = \color{blue} ,
code-for-last-col = \color{red}
}

\DeclareSymbolFontAlphabet{\mathbb}{AMSb}
\DeclareSymbolFontAlphabet{\mathbbl}{bbold}

\usepackage{aliascnt}

\newtheorem{thm}{Theorem}[section]

\newaliascnt{lem}{thm}
\newtheorem{lem}[lem]{Lemma}
\aliascntresetthe{lem}

\newaliascnt{prop}{thm}
\newtheorem{prop}[prop]{Proposition}
\aliascntresetthe{prop}

\newaliascnt{cor}{thm}
\newtheorem{cor}[cor]{Corollary}
\aliascntresetthe{cor}

\newaliascnt{conj}{thm}

\aliascntresetthe{conj}

\theoremstyle{definition}

\newaliascnt{defn}{thm}
\newtheorem{defn}[defn]{Definition}
\aliascntresetthe{defn}

\newaliascnt{nota}{thm}

\aliascntresetthe{nota}

\newaliascnt{eg}{thm}
\newtheorem{eg}[eg]{Example}
\aliascntresetthe{eg}

\newaliascnt{rem}{thm}
\newtheorem{rem}[rem]{Remark}
\aliascntresetthe{rem}

\newaliascnt{ques}{thm}

\aliascntresetthe{ques}

\theoremstyle{remark}
\newtheorem*{rem*}{Remarks}

\newcommand{\Sy}{\mathrm{Sym}}

\newcommand{\sgn}{\mathrm{sgn}}

\newcommand{\sym}[1]{\mathfrak{S}_{#1}}
\newcommand{\e}[1]{\overline{e_{#1}}}

\newcommand{\R}{\mathscr{R}}

\newcommand{\Z}{\mathbb{Z}}
\newcommand{\A}{\mathbb{A}}

\newcommand{\B}{\mathcal{B}}

\newcommand{\F}{\mathbb{F}}
\newcommand{\U}[2]{{\mathrm{U}^\wedge_{{#1}}({#2})}}

\newcommand{\res}[2]{{#1}{\downarrow_{#2}}}
\newcommand{\rad}{\mathrm{Rad}}
\DeclareMathOperator{\rank}{rank}

\numberwithin{equation}{section}
\newcommand{\rk}[2]{V^\#_{#1}({#2})}

\renewcommand{\t}{\mathfrak{t}}
\newcommand{\Sp}[2]{S^{({#1}-{#2},1^{#2})}}

\newcommand{\shortseq}[3]{0\longrightarrow #1\longrightarrow #2\longrightarrow #3\longrightarrow 0}
\begin{document}
\title[The rank varieties]{On the rank varieties of some simple modules for symmetric groups}

\author{Jialin Wang}
\address{Jialin Wang \\
School of Science \& Technology \\
Department of Mathematics \\
City St George's, University of London \\
Northampton Square \\
London EC1V 0HB \\
United Kingdom}
\email{jialin.wang@city.ac.uk}
\begin{abstract}
Let $\F$ be an algebraically closed field of odd characteristic $p$, and let $E_k$ be an elementary abelian subgroup of $\sym{kp}$ generated by $k$ disjoint $p$-cycles. Lim and the author determined the rank variety and complexity of the simple module $D^{(kp-p+1,1^{p-1})}$ when $k\not\equiv1\pmod p$ in \cite{kjw}. In this paper, we show that the result holds for every $k\geq2$. 
\end{abstract}

\maketitle

\section{Introduction}

Let $\F$ be an algebraically closed field of odd characteristic $p$ and $\sym n$ be the symmetric group of degree $n$. The simple $\F\sym n$-modules are indexed by the $p$-regular partitions of $n$. Suppose $\lambda$ is a $p$-regular partition of $n$. The simple module $D^\lambda$ is the head of the Specht module $S^\lambda$. While the Specht module $S^\lambda$ has a characteristic-free basis, given by the set of standard $\lambda$-polytabloids, little is known for the simple module.

Let $D(1):=D^{(kp-1,1)}$ and for $1\leq r\leq p-1$, let 
\[
D(r)=D^{(kp-r,1^{r})}\cong \bigwedge^rD(1).
\]
In this paper, we study these simple hook modules and their restriction to the elementary abelian subgroup $E_k$ of $\sym{kp}$ generated by $k$ disjoint $p$-cycles. Define
\begin{equation}\label{p_k}
p_k(x_1,\ldots,x_k)=
\sum_{i=1}^k\left(x_1\cdots\widehat{x_i}\cdots x_k\right)^{p-1},
\end{equation}
and for a polynomial $f\in \F[x_1,\cdots,x_k]$, let $V(f)$ be $\{\bar\alpha\in \F^k: f(\bar\alpha)=0\}$.
In \cite{kjw}, Lim and the author proved that
\[
\rk{E_k}{D(p-1)}\subseteq V(p_k),
\]
and obtained equality when $k\not\equiv1\pmod p$. Their argument does not apply when $k\equiv1\pmod p$.

We prove that, for every $k\geq2$,
\[
\rk{E_k}{D(p-1)}=V(p_k).
\]
It follows that $D(p-1)$ has complexity $k-1$. The proof does not depend on $k$ nor the dimension of $D(p-1)$. 

In Section~\ref{prelim}, we introduce notations and results used later. In Section~\ref{jordanof S}, we determine the maximal Jordan sets for hook Specht modules. In Section~\ref{Dp-1} we determine the generic Jordan types of $D(r)$ for $1\leq r \leq p-1$, and present the main result.

\section{Preliminary}\label{prelim}
Throughout, $p$ is an odd prime and $\F$ is an algebraically closed field of characteristic $p$. For the basic knowledge in this article, we refer the readers to \cite{benson2016,James78}.

\subsection{Symmetric groups and hook modules}\label{ss:sym}
Let $\sym n$ be the symmetric group permuting $\{1,
\ldots n\}$. We define $E_k$ to be the subgroup of $\sym {kp}$ generated by $p$-cycles $g_1,\ldots, g_k$ such that for $i=1,\ldots, k$,
\[g_i=((i-1)p+1,(i-1)p+2,\ldots, ip).
\]
Up to conjugacy, $E_k$ is an elementary abelian subgroup of maximal $p$-rank $k$, and
\[
N_{\sym{kp}}(E_k)/C_{\sym{kp}}(E_k)\cong \F_p^\times\wr\sym k.
\]
We say a sequence of positive integers $\lambda=(\lambda_1,\ldots, \lambda_k)$ is a partition of $n$, denoted by $\lambda\vdash n$, if $\lambda_1\geq \cdots\geq \lambda_k$ and $\lambda_1+\cdots+\lambda_k=n$. We use $l(\lambda)$ and $|\lambda|$ to denote the length and size of $\lambda$. We say $\lambda$ is $p$-regular if there are no $p$ parts of the same size and it is a hook partition if $\lambda_2=\cdots=\lambda_k=1$ or $k=1$.
 We use $|\lambda|$ and $l(\lambda)$ for the size and length of $\lambda$. A hook is a partition of the form $(n-r,1^r)$. The Young diagram of $\lambda$ is a depiction of the set $[\lambda]=\{ (i,j): 1\leq i\leq k, 1\leq j\leq \lambda_i\}$ and elements in the set are called nodes. A $\lambda$-tableau is one of the $n!$ arrays obtained by assigning each node in $[\lambda]$ by one of $1,\ldots, n$ with no repeats. We say a $\lambda$-tableau is standard if the numbers are strictly increasing in each row from left to right and in each column from top to bottom. Let $\mu=(\mu_1,\ldots,\mu_l)$ be a partition of $n$ and we say $\lambda$ dominates $\mu$, denoted as $\lambda\unrhd \mu$, if for all admissible $j$, we have 
\[\sum_{i=1}^{j} \lambda_i\geq \sum_{i=1}^j \mu_i.
\]
For each partition $\lambda$ of $n$, there is natural action of $\sym n$ on the $\lambda$-tableaux by permuting the numbers. For a $\lambda$-tableau $t$, the row (respectively column) stabilizer $R_t$ (respectively $C_t$) is the subgroup of $\sym n$ consisting of permutations which fix the rows (respectively columns) of $t$ setwise. Define an equivalence relation on the set of $\lambda$-tableaux by $t\sim s$ if $t=\sigma s$ for some $\sigma\in R_s$ and a $\lambda$-tabloid $\t$ is the equivalence class containing $t$ under the relation. The action of $\sym n$ on tabloids is the same, given by permuting the numbers.
For a $\lambda$-tableau, we define the polytabloid 
\[e_t=\sum_{\sigma\in C_t} (\sgn\sigma)\sigma \t.
\]
We say $e_t$ is standard if $t$ is standard. For a partition $\lambda\vdash n$, the Specht module $S^\lambda$ has a characteristic free basis given by the set of standard $\lambda$-polytabloids. If $\lambda$ is $p$-regular, then $S^\lambda$ has simple head $D^\lambda$.

Let $n\geq 3$ and $D(1):=D^{(n-1,1)}$ be the natural $\F\sym n$-module. For $1\leq r\leq \dim D(1)$, define the $r$-exterior power 
\[D(r):=\bigwedge^r D(1).
\]
There is a surjection $S^{(n-1,1)}\twoheadrightarrow D(1)$ and by \cite[Proposition 2.3(a)]{MZ07}, it induces a surjection $$S^{(n-r,1^r)}\cong \bigwedge^r S^{(n-1,1)}\twoheadrightarrow D(r).$$ 
 By \cite[6.3.59]{JK}, \cite{Peel}, and \cite{Danz07}, the module $D(r)$ is simple and $D(r)\cong D^{(n-r,1^r)^R}$, where $R$ denotes $p$-regularisation. In particular, when $n>p$ and $1\leq r\leq p-1$,
\[D(r)\cong D^{(n-r,1^r)}.
\]

For $i=1,\ldots,n$, let $\mathfrak t_i$ denote the
$(n-1,1)$-tabloid with $i$ in the second row. For $i>1$,
choose a tableau $t_i$ with tabloid $\mathfrak t_i$ and with
$1$ in its $(1,1)$-entry. Then
\[
e_i=e_{t_i}=\mathfrak t_i-\mathfrak t_1.
\] For $i>1$, let $e_i=e_{t_i}=\t_i-\t_1$. For any $\sigma\in \sym n$, $\sigma e_i=\t_{\sigma i}-\t_{\sigma 1}$. Then $\{e_2,\ldots,e_n\}$ is a basis of $S^{(n-1,1)}$. Suppose now $n=kp$, then $S^{(kp-1,1)}$ has a trivial socle $D^{(kp)}$ spanned by $\sum_{i=2}^{kp}e_i$, and the quotient is $D(1)$ with a basis $\{\bar e_3,\ldots,\bar e_{kp}\}$. It follows that
\[
\dim D(r)=\binom{kp-2}{r}.
\]
Notice that $p\mid \dim D(r)$ if $r\equiv p-1\pmod p$.

\subsection{Modules for $\F C_p$}\label{sectioncp}
Let $C_p$ be the cyclic group of order $p$. For $1\leq i\leq p$, let $J_i$ denote the unique indecomposable $\F C_p$-module of dimension $i$. Thus $J_1$ is the unique simple module and $J_p$ is the unique projective indecomposable module. Every finite-dimensional $\F C_p$-module $M$ has a unique decomposition
\[
M\cong\bigoplus_{i=1}^p J_i^{\oplus a_i}
\]
for some integers $a_i\geq0$. In this case, we say $M$ has Jordan type
\[
[p]^{a_p}[p-1]^{a_{p-1}}\cdots[1]^{a_1},
\]
and stable Jordan type is $[p-1]^{a_{p-1}}\cdots[1]^{a_1}$. We may also regard the Jordan type of $M$ as the partition $(p^{a_p},\ldots,1^{a_1})$ of $\dim M$.

We say two stable Jordan types are complementary if they have the form
\[
[p-1]^{a_{p-1}}\cdots[1]^{a_1}
\]
and
\[
[p-1]^{a_1}[p-2]^{a_2}\cdots[1]^{a_{p-1}}.
\]

For a partition $\lambda$, write $\bar\lambda$ for its conjugate, which is the partition with
\[
\bar\lambda_i=\bigl|\{j:\lambda_j\geq i\}\bigr|.
\]
We append zero parts when partitions of different lengths are compared.

\begin{defn}\label{D: LR seq}
Let $A=[a^0,a^1,\ldots,a^r]$ be a sequence of partitions. We say that $A$ is a Littlewood-Richardson sequence (LR sequence), if 
\begin{enumerate}[(i)]
\item $0\leq a^h_i-a^{h-1}_i\leq 1$ for $1\leq h\leq r$ and $i\geq 1$,
\item $\sum_{i\geq j}(a^{h+1}_i-a^h_i)\leq  \sum_{i\geq j}(a^h_i-a^{h-1}_i)$, for $0<h<r$ and $j\geq 1$.
\end{enumerate}
If $a^0=\lambda$, $a^r=\mu$, and $\beta$ is the partition determined by
\[
\bar\beta_h=|a^h|-|a^{h-1}|\quad(1\leq h\leq r),
\qquad \bar\beta_{r+1}=0,
\]
then $A$ is an LR sequence of type $[\lambda,\beta;\mu]$.
\end{defn}

The following theorem describes the relation of Jordan types of modules within a short exact sequence.

\begin{thm}\label{t.green}\cite[Theorem 4.1, 4.2]{klein}
Let $V_\lambda,V_\beta,V_\mu$ be $\F C_p$-modules with Jordan type $\lambda,\beta,\mu$ respectively. Then there exists a short exact sequence
\[\shortseq{V_\beta}{V_\mu}{V_\lambda}
\] of $\F C_p$-modules if and only if there exists an LR sequence of type $[\lambda,\beta;\mu]$.

\end{thm}

We recall the tableau interpretation of LR sequences. Let $\mu\subseteq\lambda$ and let $\nu$ be a partition with $|\lambda|=|\mu|+|\nu|$. A tableau of shape $\lambda\backslash\mu$ and weight $\nu$ is semistandard if its entries are weakly increasing along rows and strictly increasing down columns. The word of a tableau is obtained by reading each row from right to left, beginning with the top row and it is a lattice if the occurrences of an integer $i$ for $i>1$ are at most the occurrences of $i-1$ in any initial segment of it. In particular, a tableau is an LR tableau if this word is a lattice word.

Let $T$ be a tableau with $c$ columns and entries in $\{0,1,\ldots,r\}$. For $0\leq h\leq r$ and $1\leq a\leq b\leq c$, let $T_{a,b}(h)$ denote the number of entries equal to $h$ in columns $a,\ldots,b$, and let $T^h$ be the subtableau consisting of the entries at most $h$.

An LR sequence $A=[a^0,\ldots,a^r]$ of type $[\lambda,\beta;\mu]$ determines a tableau $T(A)$ of shape $\bar\mu$: the entries in $\{0,\ldots,h\}$ form the diagram of $\overline{a^h}$ for each $h$.

\begin{prop}\label{P: equiv LR} Let $A=[a^0,\ldots,a^r]$ be an LR sequence and $T=T(A)$. Under this correspondence, Definition~\ref{D: LR seq}(i)--(ii) becomes:
\begin{enumerate}[(i)]
\item for $1\leq h\leq r$ and $1\leq a\leq c$, $T_{a,a}(h)=0$ or $1$, i.e., the numbers are strictly increasing down the columns;
\item for each $1\leq h\leq r-1$, to the right of any column, the total number of $h$'s is at least the total number of $(h+1)$'s, i.e., for $1\leq a\leq c$, $T_{a,c}(h)\geq T_{a,c}(h+1)$.
\end{enumerate}
\end{prop}
After the zero entries are removed, conditions (i) and (ii) say precisely that the resulting tableau is semistandard and has a lattice word. Hence LR sequences of type $[\lambda,\beta;\mu]$ are in bijection with LR tableaux of shape $\bar\mu\backslash\bar\lambda$ and weight $\bar\beta$. Moreover, condition (i) says that $T^h\backslash T^{h-1}$ is a removable horizontal strip for each $1\leq h\leq r$, that is a horizontal strip whose removal leaves a Young diagram.
\begin{eg}\label{Eg: LR seq}
Let
\[
A=[a^0,a^1,a^2,a^3]
=[(2,1,1),(2,2,2,1),(3,2,2,2),(4,3,2,2)]
\]
be an LR sequence. Then $\lambda=a^0$, $\mu=a^3$, and
\[
\bar\beta
=(|a^1|-|a^0|,\ |a^2|-|a^1|,\ |a^3|-|a^2|)
=(3,2,2).
\]
Hence $\beta=(3,3,1)$, and $A$ has type
\[
[\lambda,\beta;\mu]
=[(2,1,1),(3,3,1);(4,3,2,2)].
\]
The corresponding tableau $T$ and the skew tableau $T'$ of shape
$\bar\mu\backslash\bar\lambda$ and weight $\bar\beta$ are
\[
{\ytableausetup{boxsize=1.2em,centertableaux}
T=
\begin{ytableau}
0&0&0&1\\
0&1&1&2\\
2&3\\
3
\end{ytableau}
\qquad
T'=
\begin{ytableau}
\none&\none&\none&1\\
\none&1&1&2\\
2&3\\
3
\end{ytableau}.}
\]
The tableau $T'$ is semistandard, and its word
\[
(1,2,1,1,3,2,3)
\]
is a lattice word. For example,
\[
(T_{1,1}(0),T_{1,2}(0),T_{1,3}(0),T_{1,4}(0))
=(2,3,4,4).
\]
Moreover, $T^h$ has shape $\overline{a^h}$ for $0\leq h\leq3$; in
particular, $T^3=T$.
\end{eg}

For a partition $\lambda=(\lambda_1,\ldots,\lambda_k)$ and an integer $m\geq \lambda_1$, we define $m-\lambda$ to be the partition $(m-\lambda_k,\ldots, m-\lambda_1)$.
For two partitions $\lambda$ and $\mu$, define $\lambda\sqcup\mu$ to be the partition obtained by combining their parts and arranging them in nonincreasing order.

\subsection{Complexity and rank varieties}
Let $M$ be an $\F G$-module, and let
\[\cdots\longrightarrow P_1\longrightarrow P_0\longrightarrow M\longrightarrow 0,
\]
be a minimal projective resolution. The complexity $c_G(M)$ is the smallest integer $t\geq0$ such that
\[\lim_{n\rightarrow \infty}\frac{\dim P_n}{n^t}=0.
\]
The Alperin--Evens theorem \cite{AE81} gives
\begin{thm} We have
\[c_G(M)=\max_{E\in \mathscr{E}}\{c_E(\res{M}{E})\},\]
where $\mathscr{E}$ is a set of representatives of the conjugacy classes of maximal elementary abelian $p$-subgroups of $G$.
\end{thm}
We recall Carlson's definition of the rank variety \cite{carlson}. Let $E=\langle g_1,\ldots,g_k\rangle$ be elementary abelian of rank $k$, let $J=\rad(\F E)$, and put $X_i=g_i-1$. The elements $X_1+J^2,\ldots,X_k+J^2$ form a basis of $J/J^2$. For $\alpha=(\alpha_1,\ldots,\alpha_k)\in\A^k(\F)$, write
\[
X_\alpha=\sum_{i=1}^k\alpha_iX_i,
\qquad u_\alpha=1+X_\alpha.
\]
If $\alpha\neq0$, then $u_\alpha$ has order $p$ in $(\F E)^\times$; the group $\langle u_\alpha\rangle$ is called a cyclic shifted subgroup.

For an $\F E$-module $M$, the rank variety of $M$, denoted by $\rk{E}{M}$ is
\[
\{0\}\cup
\left\{\alpha\neq0:
\res{M}{\langle u_\alpha\rangle}\text{ is not free}\right\}.
\]
Up to isomorphism, the rank variety is independent of the chosen generators of $E$ \cite[Theorem 6.5]{carlson}.
\begin{lem}[Dade's lemma {\cite[Lemma 11.8]{Dade78}}]\label{L: Dade} Let $M$ be an $\F E$-module. Then $M$ is projective if and only if $\rk{E}{M}=\{0\}$.
\end{lem}
We use the following standard properties of rank varieties.
\begin{thm}\label{T: basic rank}
Let $M$ and $N$ be $\F E$-modules. Then
\begin{enumerate}
\item $\rk{E}{M}$ is a closed homogeneous subvariety of $\rk{E}{\F}=\A^k(\F)$,
\item the dimension of $\rk{E}{M}$ is equal to the complexity $c_E(M)$,
\item $\rk{E}{M\oplus N}=\rk{E}{M}\cup\rk{E}{N}$ and $\rk{E}{M\otimes N}=\rk{E}{M}\cap\rk{E}{N}$,
\item if $M$ is indecomposable, then the image of $V^\#_E(M)$ in projective space is connected.
\end{enumerate}
\end{thm}
Let $E_k$ be the elementary $p$-subgroup of $\sym {kp}$ defined in \cref{ss:sym}. There is a natural action of $\mathrm{N}_{\sym {kp}}(E_k)/\mathrm{C}_{\sym {kp}}(E_k)$ on the rank varieties which can be described as follows.
\begin{lem}\label{l:actiononrank}
Let $M$ be an $\F\sym {kp}$-module. We have $\mathrm{N}_{\sym{kp}}(E_k)/\mathrm{C}_{\sym{kp}}(E_k)\cong \F_p^\times\wr \sym{k}$. For $\gamma\in \F_p^\times$ in the $i$th component, $\sigma\in\sym{k}$ and $\alpha\in \rk{E_k}{M}$,
\begin{align*}
  \gamma\cdot \alpha&=(\alpha_1,\ldots,\gamma\alpha_i,\ldots,\alpha_k),\\
  \sigma\cdot \alpha&=(\alpha_{\sigma^{-1}(1)},\alpha_{\sigma^{-1}(2)},\ldots,\alpha_{\sigma^{-1}(k)}).
\end{align*}
\end{lem}

\subsection{Generic and maximal Jordan types}
Let $E=\langle g_1,\ldots,g_k\rangle$. Followed the notation above, for an $\F E$-module $M$ and an element $x\in J$, the matrix representation of $x$ with respect to some basis of $M$ is a sum of Jordan blocks. We say $x$ on $M$ has Jordan type
\[
[p]^{a_p}\cdots[1]^{a_1}
\]
if it has exactly $a_i$ Jordan blocks of size $i$. 
We say $x\in J$ has maximal Jordan type if the Jordan type of $x$ is maximal with respect to the dominance order among all the Jordan types of elements in $J$. The maximal Jordan set $\U{E}{M}$ is the set of elements in $\F E$ of maximal Jordan type. In particular, we have $\U{E}{M}\subseteq (J/J^2)\backslash\{0\}$ and we may write $x\in \U{E}{M}$ as $x=\sum_i \alpha_iX_i+J^2$ for some $(\alpha_1,\ldots, \alpha
_k)\in \A^k(\F)$. Then, we may regard $\U{E}{M}\subseteq \A^k(\F)$.
The Jordan type attained on a dense open subset of $\A^k(\F)$ (or roughly speaking almost every point in $\A^k(\F)$) is called the generic Jordan type of $M$. In particular, we have the following theorem.
  
\begin{thm}\cite[\S 4]{FPS}\label{t:genericmax} 
Let $M$ be an $\F E$-module. 
All elements of maximal Jordan type on $M$ have the same Jordan type and the type is the same as the generic Jordan type of $M$.
\end{thm} 
An $\F E$-module is generically free if its generic Jordan type is $[p]^d$ for some $d$. This can occur only when $p$ divides $\dim M$. If $M$ is generically free, then $\res{M}{\langle u_\alpha\rangle}$ is free precisely when $\alpha\in\U{E}{M}$. For $U\subseteq\A^k(\F)$, write $U^c$ for its complement. It follows that
\begin{cor}\label{lem:complement}
For an $\F E$ module $M$, if $M$ is generically free, then $\U{E}{M}=\rk{E}{M}^c$.
\end{cor}
Stable and complementary stable Jordan types of $\F C_p$-modules are defined as in \cref{sectioncp}. For $\F E$-modules $M_1, M_2, M_3$, with the exact sequence \[\shortseq{M_1}{M_2}{M_3},
\]
if $M_2$ is generically free, then $M_1$ and $M_3$ have complementary stable Jordan types. If $M_1$ or $M_3$ is generically free, then the other two modules have the same stable generic Jordan types.
\begin{lem}\label{u&seq}
Suppose we have the following short exact sequence for $\F E$-modules $M_1, M_2, M_3$:
\[\shortseq{M_1}{M_2}{M_3}.
\]
Let $\{i_1,i_2,i_3\}=\{1,2,3\}$. If $M_{i_1}$ is generically free, then we have 
\[\U{E}{M_{i_1}}\cap\U{E}{M_{i_2}}=\U{E}{M_{i_1}}\cap\U{E}{M_{i_3}}.\]
\end{lem}

\begin{proof}
Fix $\alpha\in \U{E}{M_{i_1}}$. Since $M_{i_1}$ is generically free, its restriction to $\langle u_\alpha\rangle$ is free. The algebra $\F\langle u_\alpha\rangle\cong \F C_p$ is self-injective. Hence, if $i_1=1$ or $3$, the restricted exact sequence splits and the other two terms have the same stable Jordan type. If $i_1=2$, the restrictions of $M_1$ and $M_3$ have complementary stable Jordan types. The same relation holds for their generic Jordan types. Thus, on $\U{E}{M_{i_1}}$, one of the other two modules has maximal Jordan type if and only if the other does. 
\end{proof}

For $1\leq r\leq p-1$, the $r$-th exterior power is the Schur functor corresponding to $(1^r)$. The following result is \cite[Corollary 4.6.2]{benson2016}.
\begin{prop}\label{p:wedgeJT}
For an $\F E$-module $M$ and an integer $1\leq r\leq p-1$, the generic Jordan type of $\bigwedge^r(M)$ is the $r$-th exterior power of the generic Jordan type of $M$ and $\U{E}{\bigwedge^r(M)}\supseteq \U{E}{M}$.
\end{prop}

\section{Maximal Jordan sets of hook Specht modules}\label{jordanof S}
For $1\leq i\leq k$, recall that
\[
g_i=((i-1)p+1,\ldots,ip)\in\sym{kp},
\qquad X_i=g_i-1.
\]
We determine the maximal Jordan sets of the hook Specht modules $S^{(n-r,1^r)}$, when $n\in\{kp,kp+1\}$ and $1\leq r\leq p-1$. By \cite{kj09}, the stable generic Jordan type of $\res{S^{(n-r,1^r)}}{E_k}$ is
\[
\begin{cases}
[1],&n=kp\text{ and }r\text{ is even},\\
[p-1],&n=kp\text{ and }r\text{ is odd},\\
\emptyset,&n=kp+1.
\end{cases}
\]

For $2\leq i\leq k$ and $1\leq s\leq p$, write $e_{i,s}=e_{(i-1)p+s}$. In $S^{(kp-1,1)}$, for $2\leq i\leq k$, define
\[
b_1=e_2,\qquad b_i=e_{i,1}-e_2\quad,
\]
and let
\[
\B_1=\{b_1,X_1b_1,\ldots,X_1^{p-2}b_1\},\qquad
\B_i=\{b_i,X_ib_i,\ldots,X_i^{p-1}b_i\}.
\]
Lastly let $\B=\bigcup_{i=1}^k \B_i$.

For $S^{(kp,1)}$, we define $\mathfrak{t}_{i,r}'$, $ e_{i,r}'$, $b_1'$, $b_i'$, $\B_i'$ and $\B'$ accordingly.

\begin{lem}\label{l.basis}
For $1\leq r\leq p-2$, $1\leq m\leq p-1$, and $2\leq i\leq k$,
\[
X_1^rb_1=\sum_{s=1}^{r+1}(-1)^{r-s+1}{r+1\choose s}e_{s+1},\qquad
X_i^mb_i=\sum_{s=0}^{m}(-1)^{m-s}{m\choose s}e_{i,s+1}.
\]
In particular, $\B$ is a basis of $S^{(kp-1,1)}$.
\end{lem}

\begin{proof}
We prove the first identity by induction on $r$. When $r=1$, we have
$$
X_1b_1=g_1e_2-e_2=e_3-2e_2.
$$
Suppose the first identity holds for some $1\leq r\leq p-3$. We have:
\begin{align*}
X_1(X_1^rb_1)
&=X_1(\sum_{s=1}^{r+1} (-1)^{r-s+1} {r+1\choose s} e_{s+1})\\
&=\sum_{s=1}^{r+1} (-1)^{r-s+1} {r+1\choose s} (e_{s+2}-e_2)-(\sum_{s=1}^{r+1} (-1)^{r-s+1} {r+1\choose s} e_{s+1})\\
&=(-1)^{r+1}(r+2)e_2+\sum_{s=1}^{r} (-1)^{r-s+1} {r+2\choose s+1} e_{s+2}+e_{r+3}\\
&=\sum_{s=1}^{r+2} (-1)^{r+1-s+1} {r+2\choose s} e_{s+1},
\end{align*}
as desired. The proof for the second formula is similar and thus omitted here.
Suppose that $T$ is the matrix such that 
\[(b_1,\ldots,X_1^{p-2}b_1,b_2,\ldots, X_2^{p-1}b_2,\ldots, X_k^{p-1}b_k)^{\top}=T(e_2,\ldots,e_{kp})^\top.
\]
Then $T$ is a lower triangular matrix with nonzero diagonal entries. So $\B$ is linearly independent with $|\B|=kp-1$, it follows that $\B$ is a basis.
\end{proof}
 
We state the same result for $S^{(kp,1)}$. Its proof is similar to that of Lemma \ref{l.basis} and thus omitted.

\begin{lem}\label{l.basis'}
For $1\leq r\leq p-2$, $1\leq m\leq p-1$, and $2\leq i\leq k$,
\[
X_1^rb'_1
=\sum_{s=1}^{r+1}(-1)^{r-s+1}{r+1\choose s}e'_{s+1},\qquad
X_i^mb'_i
=\sum_{s=0}^{m}(-1)^{m-s}{m\choose s}e'_{i,s+1}.
\]
Moreover, $\B'$ is a basis of $S^{(kp,1)}$.
\end{lem}

\begin{lem}\label{l.action1}
Let $2\leq i\leq k$, $0\leq r\leq p-1$, and $1\leq j\leq k$. Then
\[
X_j(X_i^rb_i)=
\begin{cases}
X_i^{r+1}b_i,&j=i\text{ and }0\leq r\leq p-2,\\
-b_1-X_1b_1,&j=1\text{ and }r=0,\\
0,&\text{otherwise}.
\end{cases}
\]
We also have $
X_1(X_1^{p-2}b_1)=0$,
and if $j\neq 1$, $X_j(X_1^rb_1)=0.$
Consequently, $X_iX_j$ annihilates $S^{(kp-1,1)}$ whenever $i\ne j$.
\end{lem}

\begin{proof}
Suppose $i$ and $j$ are distinct elements of $\{2,\ldots,k\}$, then $g_j$ fixes every $e_{i,s}$. For $j=1$, one has $X_1e_{i,s}=-e_2$. Hence $X_1b_i=-b_1-X_1b_1$, whereas the binomial formula in Lemma \ref{l.basis} gives $X_1X_i^rb_i=0$ for $r\geq1$. For $r=p-2$, we have $X_1^{p-2}b_1=\sum_{s=2}^p e_s$ which is fixed under $g_1$.
The remaining identities follow directly from
\[
X_i^p=(g_i-1)^p=g_i^p-1=0.
\]
\end{proof}
A similar result holds for $S^{(kp,1)}$, of which the proof is omitted.
\begin{lem}\label{l.action2}
Let $2\leq i\leq k$, $0\leq r\leq p-1$, and $1\leq j\leq k$. Then
\[
X_j(X_i^rb_i')=
\begin{cases}
X_i^{r+1}b_i',&j=i\text{ and }0\leq r\leq p-2,\\
-b_1'-X_1b_1',&j=1\text{ and }r=0,\\
0,&\text{otherwise}.
\end{cases}
\]
We also have $
X_1(X_1^{p-2}b_1')=0$,
and if $j\neq 1$, $X_j(X_1^rb_1')=0.$
Consequently, $X_iX_j$ annihilates $S^{(kp,1)}$ whenever $i\ne j$.
\end{lem}

Fix $0\ne\alpha=(\alpha_1,\ldots,\alpha_k)\in\A^k(\F)$ and write $X_\alpha=\sum_{i=1}^k\alpha_iX_i$. For $d\geq1$ and $a\in \F^\times$, let $N_d(a)$ denote the $d\times d$ matrix with $a$ on the subdiagonal and zero elsewhere. Let
\[
c(a)=(-a,-a,0,\ldots,0)\in\F^{p-1},
\]
and define
\[
C_d(a)=c(a)^{\mathsf T}(1,0,\ldots,0)\in M_{p-1,d}(\F).
\]
Representing $X_\alpha$ with respect to $\B$ and $\B'$, we have
\[
L=[X_\alpha]_{\B}=
\begin{pmatrix}
N_{p-1}(\alpha_1)&C_p(\alpha_1)&\cdots&C_p(\alpha_1)\\
&N_p(\alpha_2)&&\\
\vdots&&\ddots&\\
&&&N_p(\alpha_k)
\end{pmatrix},
\]
and
\[
L'=[X_\alpha]_{\B'}=
\begin{pmatrix}
N_{p-1}(\alpha_1)&C_p(\alpha_1)&\cdots&C_p(\alpha_1)&C_1(\alpha_1)\\
&N_p(\alpha_2)&&&\\
\vdots&&\ddots&&\vdots\\
&&&N_p(\alpha_k)&\\
0&0&\cdots&0&0
\end{pmatrix}.
\]
All nonzero off-diagonal blocks lie in the first block row.

Define
\[
f_k=x_1\cdots x_k\in\F[x_1,\ldots,x_k].
\]
Then
\[
V(f_k)^c=\{\alpha\in\A^k(\F):\alpha_i\ne0\text{ for all }i\}.
\]
\begin{lem}\label{l.rank}
For $k\geq2$,
\begin{enumerate}[(i)]
\item $\rank L\leq(k-1)(p-1)+p-2$, with equality if and only if $f_k(\alpha)\ne0$;
\item $\rank L'\leq k(p-1)$, with equality if and only if $f_k(\alpha)\ne0$.
\end{enumerate}
\end{lem}

\begin{proof}
If $\alpha_1=0$, then the first diagonal block and all the blocks $C_d(\alpha_1)$ vanish, so neither upper bound is attained. Suppose that $\alpha_1\ne0$, and let
\[
z=\bigl|\{i:2\leq i\leq k\text{ and }\alpha_i=0\}\bigr|.
\]
The image of $N_{p-1}(\alpha_1)$ has dimension $p-2$, and $c(\alpha_1)$ does not belong to this image. In a block column for which $\alpha_i\ne0$, the vector $c(\alpha_1)$ is tied to a nonzero column of $N_p(\alpha_i)$ and does not increase the rank of that block column. If $\alpha_i=0$, the same block column contributes the vector $c(\alpha_1)$ alone. All such vectors are equal, so together they contribute exactly one additional dimension when $z>0$. It follows that
\[
\rank L=
\begin{cases}
p-2+(k-1)(p-1),&z=0,\\
p-2+(k-1-z)(p-1)+1,&z>0,
\end{cases}
\]
For $L'$, the final block column already contributes $c(\alpha_1)$ independently of the diagonal blocks, and hence
\[
\rank L'=p-2+(k-1-z)(p-1)+1.
\]
Since $p$ is odd, the displayed ranks attain the stated upper bounds exactly when $z=0$. Together with the case $\alpha_1=0$, this proves both assertions.
\end{proof}

The generic Jordan types of $\res{S^{(kp-1,1)}}{E_k}$ and $\res{S^{(kp,1)}}{E_k}$ are $[p]^{k-1}[p-1]$ and $[p]^k$, respectively.

\begin{thm}\label{t.umaxS}
We have
\[\U{E_k}{S^{(kp-1,1)}}=\U{E_k}{S^{(kp,1)}}=V(f_k)^c.\]
In particular, $\rk{E_k}{S^{(kp,1)}}=V(f_k)$.
\end{thm}

\begin{proof}
The restriction of $S^{(kp,1)}$ to $E_k$ is generically free. Consequently, $\alpha$ belongs to its maximal Jordan set if and only if
\[
\rank L'=k(p-1).
\]
Part (ii) of \cref{l.rank} therefore gives
\[
\U{E_k}{S^{(kp,1)}}=V(f_k)^c
\qquad\text{and}\qquad
\rk{E_k}{S^{(kp,1)}}=V(f_k).
\]

Now consider $S^{(kp-1,1)}$. Since $X_\alpha^p=0$, all Jordan blocks of $X_\alpha$ have size at most $p$. An operator on a space of dimension $kp-1$ having rank $(k-1)(p-1)+p-2$ has exactly $k$ Jordan blocks. Their dimensions sum to $kp-1$, so their sizes must be $p,\ldots,p,p-1$. Thus $X_\alpha$ has maximal Jordan type precisely when $\rank L=(k-1)(p-1)+p-2$. Part (i) of \cref{l.rank} now yields
\[
\U{E_k}{S^{(kp-1,1)}}=V(f_k)^c.
\]
\end{proof}

For $1\leq r\leq p-1$, we have $S^{(kp+1-r,1^r)}\cong \bigwedge\nolimits ^r S^{(kp,1)}$. Proposition~\ref{p:wedgeJT} and Theorem~\ref{t.umaxS} therefore give the following corollary.
\begin{cor}\label{freeskp}
For $1\leq r\leq p-1$, we have $\U{E_k}{S^{(kp+1-r,1^r)}}\supseteq V(f_k)^c$ and 
$$\rk{E_k}{S^{(kp+1-r,1^r)}}\subseteq V(f_k).$$
\end{cor}
An analogous statement holds for the hook Specht modules of $\sym{kp}$.
\begin{thm}\label{t.S(r)}
For $1\leq r\leq p-1$, we have $\U{E_k}{\Sp{kp}{r}}\supseteq V(f_k)^c$.
\end{thm}

\begin{proof}
We prove the theorem by induction on $r$. When $r=1$, it holds by \cref{t.umaxS}. For $2\leq r \leq p-1$, we have the following short exact sequence by branching rule:
\[0\longrightarrow \Sp{kp}{r}\longrightarrow \res{S^{(kp+1-r,1^r)}}{\sym{kp}}\longrightarrow S^{(kp-r+1,1^{r-1})}\longrightarrow 0.
\]
Since $\res{(\res{\Sp{kp+1}{r}}{\sym{kp}})}{E_k} \cong\res{S^{(kp+1-r,1^r)}}{E_k}$, we have
\[0\longrightarrow \res{\Sp{kp}{r}}{E_k} \longrightarrow \res{S^{(kp+1-r,1^r)}}{E_k}\longrightarrow\res{S^{(kp-r+1,1^{r-1})}}{E_k} \longrightarrow 0.
\]
The module $\res{S^{(kp+1-r,1^r)}}{E_k}$ is generically free, and \cref{freeskp} gives $\U{E_k}{S^{(kp+1-r,1^r)}}\supseteq V(f_k)^c$. The induction hypothesis gives $\U{E_k}{S^{(kp-r+1,1^{r-1})}}\supseteq V(f_k)^c$. Applying \cref{u&seq}, we obtain
\[V(f_k)^c\subseteq \U{E_k}{S^{(kp+1-r,1^r)}}\cap\U{E_k}{S^{(kp-r+1,1^{r-1})}}=\U{E_k}{S^{(kp+1-r,1^r)}}\cap\U{E_k}{\Sp{kp}{r}}.
\]
Therefore, \[\U{E_k}{\Sp{kp}{r}}\supseteq V(f_k)^c.\]
\end{proof}

\section{The simple modules $D(p-1)$}\label{Dp-1}
Fix $k\geq2$. For $1\leq r\leq p-1$, write
\[
D(r)=\bigwedge^r D^{(kp-1,1)}\cong D^{(kp-r,1^r)}.
\]
We retain the polynomial $p_k$ defined in \eqref{p_k}.
By \cite[Lemma 2.2]{kjw}, when $k\geq 3$, the polynomial $p_k$ is irreducible and the variety $V(p_k)$ is irreducible of dimension $k-1$.
Let $B$ be the basis of $D(1)$ used in \cite{kjw}, and let $S=[X_\alpha]_B$.

\begin{lem}\label{l:rank}\cite[Lemma 4.4]{kjw}
Let $k\geq 2$ and $S=[X_\alpha]_B$. Then
\begin{enumerate}[(i)]
\item $\rank(S)\leq (k-1)(p-1)+p-3$, with equality if and only if all coordinates $\alpha_i$ are nonzero;
\item if $p\geq 5$ and all coordinates $\alpha_i$ are nonzero, then $\rank(S^{p-3})=3k-2$ and $\rank(S^{p-2})=2k-2$;
\item if all coordinates $\alpha_i$ are nonzero, then $\rank(S^{p-1})\leq k-1$, with equality if and only if $p_k(\alpha_1,\alpha_2,\ldots,\alpha_k)\neq 0$.
\end{enumerate}
\end{lem}
\begin{rem}
We record two corrections to \cite{kjw}. The formula in \cite[Lemma 4.2(ii)]{kjw} should read
$$X_1^{p-3}b_1=\sum_{s=1}^{p-2}s\e {s+2}.
$$
In the proof of \cite[Lemma 4.4]{kjw}, if $\alpha_i=0$ for some $2\leq i\leq k$, the correct upper bound is $\rank(S)\leq(k-2)(p-1)+p-1$, rather than $(k-2)(p-1)+p-2$.
\end{rem}
Recall that $f_k=x_1\cdots x_k$. The next theorem follows from the lemma above and Theorem 4.5 in \cite{kjw}.

\begin{thm}\label{t.JTD1}
The generic Jordan type of $\res{D(1)}{E_k}$ is $[p]^{k-1}[p-2]$ and we have $$\U{E_k}{D(1)}=(V(f_k)\cup V(p_k))^c.$$ Furthermore, if $k\geq 2$, and $\alpha \in {V(f_k)}^c\cap V(p_k)$, then $\res{D(1)}{\langle u_\alpha\rangle}$ has Jordan type $[p]^{k-2}[p-1]^{2}$.
\end{thm}
For a finite group $G$ and an $\F G$-module $M$, let $\Sy^s(M)$ and $\bigwedge^s(M)$ denote the $s$th symmetric and exterior powers of $M$, respectively. We use the following results from \cite{benson2016}.
\begin{thm}\cite{benson2016}\label{propsym}
For $1\leq i\leq p$, let $J_i$ denote the indecomposable $\F C_p$-module of dimension $i$. We have
\begin{enumerate}
\item if $s<p$ and $s+i>p$, then $\Sy^s(J_i)$ is projective, and hence isomorphic to a direct sum of copies of $J_p$;
\item if $s+i\leq p+1$, then $\Sy^s(J_i)\cong \bigwedge^{s}(J_{s+i-1})$;
\item for $1\leq i\leq p$, $J_i\otimes J_p\cong J_p^{\oplus i}$.
\end{enumerate}
\end{thm}
\begin{lem}\label{l.Jp-2}
\begin{enumerate}[(i)]
\item For $1\leq r\leq p-1$, we have
\[\bigwedge^r(J_{p-1})\cong \begin{cases}
 J_{1}\oplus \text{copies of }J_p,&r\text{ is even},\\
 J_{p-1}\oplus \text{copies of }J_p,&r\text{ is odd}.
\end{cases}
\]
\item For $1\leq r\leq p-2$, we have
\[\bigwedge^r(J_{p-2})\cong \begin{cases}
 J_{r+1}\oplus \text{copies of }J_p,&r\text{ is even},\\
 J_{p-r-1}\oplus \text{copies of }J_p,&r\text{ is odd}.
\end{cases}
\]
\end{enumerate}
\end{lem}
\begin{proof}
For part (i), by \cref{propsym}, $\bigwedge^r(J_{p-1})\cong \Sy^r(J_{p-r})$.
Consider the short exact sequence \[\shortseq{J_{p-r}}{J_{p-r+1}}{J_1}.\] The exact sequences of symmetric and exterior powers in \cite[Lemma 3]{bensonlim} give
\[\shortseq{\Sy^r(J_{p-r})}{\Sy^r(J_{p-r+1})}{\Sy^{r-1}(J_{p-r+1})\otimes \bigwedge^1(J_1)},
\]
since the remaining terms on the right vanish, as $\bigwedge^i(J_1)=0$ for $i>1$. By \cref{propsym} again, $\Sy^r(J_{p-r+1})$ is projective. It follows that $\Sy^r(J_{p-r})$ and $\Sy^{r-1}(J_{p-r+1})$ have complementary Jordan types. Since $\Sy^1(J_{p-1})\cong J_{p-1}$, we have 
\[\bigwedge\nolimits ^2(J_{p-1})\cong\Sy^2(J_{p-2})\cong J_1\oplus \text{copies of $J_p$},\] 
and induction on $r$ proves (i).

For (ii), work in the representation ring $\R=a(\F C_p)$. Following \cite[Section 2.7]{benson2016}, let $\hat{\R}=\R[q]/(q^2-J_2q+1)$. In this ring, $J_n$ corresponds to $q^{n-1}+q^{n-3}+\cdots+q^{-n+3}+q^{-n+1}$ for $1\leq n\leq p$. Define the Gaussian polynomial
\[g_{n,r}(q)=\frac{(q^n-q^{-n})(q^{n-1}-q^{-n+1})\cdots (q^{n-r+1}-q^{-n+r-1})}{(q^r-q^{-r})(q^{r-1}-q^{-r+1})\cdots (q-q^{-1})}.\]
For $r<p$, we have in $\hat\R$
\[\bigwedge^r{J_n}=g_{n,r}(q),
\]
and
\begin{equation}\label{gaussian}
g_{n,r}(q)=q^{n-r}g_{n-1,r-1}(q)+q^{-r}g_{n-1,r}(q).
\end{equation}
We prove (ii) by induction on $r$. The case $r=1$ is immediate. Since only the stable Jordan type is relevant, we work in the quotient
\[\bar\R=\Z[q,q^{-1}]/(q-1)(q^{p-1}+q^{p-3}+\cdots+q^{-p+1}).
\]
The representation ring $\R$ can be identified with  the subring of $\bar\R$ consisting of fixed points of the automorphism interchanging $q$ and $q^{-1}$, and $J_1,\ldots,J_p$ remain linearly independent. In $\bar\R$, we have $q^p=q^{-p}$, and the elements $q^{p}+q^{p-2}+\cdots+q^{-p+2}$, $q^{p-1}+q^{p-3}+\cdots+q^{-p+1}$, and $q^{p-2}+q^{p-4}+\cdots+q^{-p}$ all correspond to $J_p$. We use $j(q)$ to denote any of these representatives of $J_p$ in $\bar\R$; then $qj(q)=j(q)$. We consider the two parities of $r$.

If $r$ is even, then by induction hypothesis, $\bigwedge^{r-1}J_{p-2}\cong J_{p-r}\oplus a J_p$ and $\bigwedge^r J_{p-1}\cong J_1\oplus b J_p$ for some $a,b$, and they correspond to $q^{p-r-1}+q^{p-r-3}+\cdots+q^{-p+r+1}+a j(q)$ and $1+bj(q)$ in $\bar\R$ respectively. By \cref{gaussian}, we have
\begin{align*}
g_{p-2,r}(q)=&q^rg_{p-1,r}(q)-q^{p-1}g_{p-2,r-1}\\
                  =&q^r+q^rbj(q)-q^{p-1}(q^{p-r-1}+q^{p-r-3}+\cdots+q^{-p+r+1}+a j(q))\\
                  =&q^r-(q^{2p-r-2}+q^{2p-r-4}+\cdots+q^{r})+(b-a)j(q)\\
                  =&-(q^{2p-r-2}+q^{2p-r-4}+\cdots+q^{p+1}+q^{p-1}+\cdots+q^{r+2})+(b-a)j(q)\\
                  =&(b-a-1)j(q)+(q^{p-1}+\cdots+q^{-p+1})\\
                    &-(q^{p-1}+\cdots+q^{r+2}+q^{-r-2}+q^{-r-4}+\cdots+q^{-p+1})\\
                  =&q^{r}+q^{r-2}+\cdots+q^{-r}+(b-a-1)j(q),
\end{align*}
which corresponds to $J_{r+1}\oplus (b-a-1)J_p$.

 If $r$ is odd, then by induction hypothesis again $\bigwedge^{r-1}J_{p-2}\cong J_{r}\oplus a J_p$ and $\bigwedge^r J_{p-1}\cong J_{p-1}\oplus b J_p$ for some $a,b$, and they correspond to $q^{r-1}+q^{r-3}+\cdots+q^{-r+1}+a j(q)$ and $q^{p-2}+\cdots+q^{-p+2}+bj(q)$ in $\bar\R$ respectively. By \cref{gaussian}, we have
\begin{align*}
g_{p-2,r}(q)=&q^rg_{p-1,r}(q)-q^{p-1}g_{p-2,r-1}\\
                  =&q^r(q^{p-2}+\cdots+q^{-p+2}+bj(q))-q^{p-1}(q^{r-1}+q^{r-3}+\cdots+q^{-r+1}+a j(q))\\
                  =&q^{p+r-2}+\cdots+q^{-p+r+2}+bj(q)-(q^{p+r-2}+q^{p+r-4}+\cdots+q^{p-r}+a j(q))\\
                  =&q^{p-r-2}+q^{p-r-4}+\cdots+q^{-p+r+2}+(b-a)j(q),
\end{align*}
which corresponds to $J_{p-r-1}\oplus (b-a)J_p$.
\end{proof}
Since an $\F C_p$-module is determined up to isomorphism by its Jordan type, direct sums, tensor products, and symmetric and exterior powers may be applied to Jordan types. 

Next we determine the generic Jordan types of the modules $D(r)$ when $1\leq r\leq p-1$.

\begin{lem}\label{l.gJTDr}
\begin{enumerate}[(i)]
\item For $1\leq r\leq p-2$, the stable generic Jordan type of $\res{D(r)}{E_k}$ is $[r+1]$ if $r$ is even, and $[p-1-r]$ if $r$ is odd. 
\item For $k\geq 2$ and $r=p-1$, the generic Jordan type of $\res{D(p-1)}{E_k}$ is $[p]^m$ where $mp=\binom{kp-2}{p-1}$.
\end{enumerate}
\end{lem}

\begin{proof}
Suppose first that $1\leq r\leq p-2$. By Proposition~\ref{p:wedgeJT} and Theorem~\ref{t.JTD1}, the generic Jordan type of $\res{D(r)}{E_k}$ is the Jordan type of
\[
\bigwedge\nolimits ^r(J_p^{\oplus k-1}\oplus J_{p-2}).
\]
By the following decomposition (see \cite[Lemma 1.14.3]{benson2016}), we have
\[
\bigwedge\nolimits^r
   \bigl(J_p^{\oplus k-1}\oplus J_{p-2}\bigr)
\cong
\bigoplus_{i=0}^r
   \bigwedge\nolimits^i J_{p-2}
   \otimes
   \bigwedge\nolimits^{r-i}
      \bigl(J_p^{\oplus k-1}\bigr).
\]
For $i<r$, the second tensor factor is projective, and hence so is the corresponding summand. The stable generic Jordan type is therefore the stable Jordan type of $\bigwedge^rJ_{p-2}$, which is given by Lemma \ref{l.Jp-2}(ii).

For $r=p-1$, consider
\[
N=\bigwedge\nolimits ^{p-1}(J_{p-2}\oplus J_p^{\oplus k-1})
\cong\bigoplus_{i=0}^{p-1}
\bigwedge\nolimits ^iJ_{p-2}\otimes
\bigwedge\nolimits ^{p-1-i}(J_p^{\oplus k-1}).
\]
Every summand with $i<p-1$ is projective, while the summand with $i=p-1$ is zero because $\dim J_{p-2}=p-2$. Hence $N$ is projective. Since $\dim D(p-1)=\binom{kp-2}{p-1}$, its generic Jordan type is $[p]^m$ with $mp=\binom{kp-2}{p-1}$.
\end{proof}

For a partition $\lambda=(\lambda_1,\ldots,\lambda_t)$ and an integer $m\geq\lambda_1$, recall that
\[
m-\lambda=(m-\lambda_t,\ldots,m-\lambda_1).
\]
\begin{lem}\label{l.lrprop2}
Let $A=[a^0,\ldots,a^r]$ be an LR sequence of type $[\lambda,\beta;\mu]$, and let $2\leq m\leq p-1$.
\begin{enumerate}[(i)]
\item If $\mu=(p^b,1)$ and $\beta=p-\beta'$ for some partition $\beta'\vdash m$ with $\beta'\neq(m)$, then
\[
\lambda=(p^h)\sqcup\delta
\]
for some partition $\delta\vdash m+1$ with $\delta\neq(m+1)$, where $h=b-\ell(\delta)$.
\item If $\mu=(p^b,p-1)$ and $\beta\vdash m$ with $\beta\neq(m)$, then
\[
\lambda=(p^h)\sqcup(p-\delta)
\]
for some partition $\delta\vdash m+1$ with $\delta\neq(m+1)$, where $h=b-\ell(\delta)+1$.
\end{enumerate}
\end{lem}
The proof is given in Appendix A.

\begin{thm}\label{t.umaxDr}
Suppose that $k\geq 2$. For $1\leq r\leq p-1$, we have $\U{E_k}{D(r)}^c\supseteq V(f_k)^c\cap V(p_k)$.
\end{thm}

\begin{proof}
Let $\alpha\in V(f_k)^c\cap V(p_k)$ and $1\leq r\leq p-1$. Define 
\[Q_r:=\left\{\begin{array}{ll}
\{\mu: \mu\vdash r+1, \mu\neq (r+1)\}, &\text{ if $r$ is even,} \\
\{p-\mu: \mu\vdash r+1, \mu\neq (r+1)\}, &\text{ if $r$ is odd.} \end{array}\right.
       \]
For $1\leq r\leq p-1$, the maximal stable Jordan type of $\res{D(r)}{E_k}$ does not lie in $Q_r$. Therefore, it is enough to prove that the stable Jordan type of $\res{D(r)}{\langle u_\alpha\rangle}$ belongs to $Q_r$. We prove this by induction on $r$.

The case $r=1$ follows from \cref{t.JTD1}. Assume the assertion for $r-1$. The composition factors of $S^{(kp-r,1^r)}$ are $D(r-1)$ and $D(r)$, with head $D(r)$. Restriction to $\langle u_\alpha\rangle$ gives the short exact sequence
\[\shortseq{\res{D(r-1)}{\langle u_\alpha\rangle}}{\res{S^{(kp-r,1^r)}}{\langle u_\alpha\rangle}}{\res{D(r)}{\langle u_\alpha\rangle}}.
\] By \cref{t.S(r)}, $\alpha\in \U{E_k}{S^{(kp-r,1^r)}}$ and  the stable Jordan type of $\res{S^{(kp-r,1^r)}}{\langle u_\alpha\rangle}$ is $[1]$ if $r$ is even and $[p-1]$ otherwise. 
Let $\beta\in Q_{r-1}$ be the stable Jordan type of $\res{D(r-1)}{\langle u_\alpha\rangle}$, and let $\lambda$ be the stable Jordan type of $\res{D(r)}{\langle u_\alpha\rangle}$. There are $\F\langle u_\alpha\rangle$-submodules $X$ and $P$ such that $X$ has Jordan type $\beta$, $P$ is projective, and
\[\res{D(r-1)}{\langle u_\alpha\rangle}\cong X\oplus P.
\]
Let $\phi$ denote the injection into $\res{S^{(kp-r,1^r)}}{\langle u_\alpha\rangle}$. Since $\F\langle u_\alpha\rangle\cong \F C_p$ is self-injective and $\phi(P)$ is projective, we may choose a complement $Y$ to $\phi(P)$ containing $\phi(X)$. The stable Jordan type of $Y$ is the same as that of the restricted Specht module, and there is a short exact sequence
\[\shortseq{X}{Y}{\res{D(r)}{\langle u_\alpha\rangle}}.
\]
If $r$ is even, \cref{t.green} gives an LR sequence of type $[(p^a)\sqcup\lambda,\beta;(p^b,1)]$ for some nonnegative integers $a$ and $b$. Hence Lemma \ref{l.lrprop2} gives $\lambda\in Q_r$. If $r$ is odd, the same argument uses an LR sequence of type $[(p^b)\sqcup\lambda,\beta;(p^a,p-1)]$ and again gives $\lambda\in Q_r$.
\end{proof}

\begin{lem}\label{l.inclusion}
Suppose $k\geq 2$. We have
$$\rk{E_k}{D(p-1)}\cap V(f_k)=\bigcup_{1\leq i\neq j\leq k}V(x_i,x_j).$$
\end{lem}

\begin{proof}
Since $V(f_k)=\bigcup_{1\leq i\leq k}V(x_i)$, it suffices to prove that
\[\rk{E_k}{D(p-1)}\cap V(x_i)=\bigcup_{1\leq j\leq k,j\neq i}V(x_i,x_j).\]
The proof of \cite[Lemma 5.2]{kjw} gives
\[\rk{E_k}{D(p-1)}\cap V(x_i)\subseteq\bigcup_{1\leq j\leq k,j\neq i}V(x_i,x_j).\] 
If $k=2$, then for $i=1,2$,
\[\rk{E_k}{D(p-1)}\cap V(x_i)\subseteq V(p_k)\cap V(x_i)=V(x_1,x_2)=\{0\}.
\]
So the equality holds. Now assume $k>2$. It remains to show that if two coordinates of $\alpha$ vanish, then $\res{D(p-1)}{\langle u_\alpha\rangle}$ is not free. By Lemma \ref{l:actiononrank}, we may assume that $\alpha_{k-1}=\alpha_k=0$. Let $\alpha'=(\alpha_1,\ldots,\alpha_{k-2})$. Then $u_\alpha=u_{\alpha'}\in\F E_{k-2}$ and
\[\res{(\res{D(p-1)}{E_{k-2}})}{\langle u_{\alpha'}\rangle}\cong \res{D(p-1)}{\langle u_{\alpha'}\rangle}= \res{D(p-1)}{\langle u_\alpha\rangle}.
\]
Moreover,
\[
\bigwedge^{p-1}(\res{D(1)}{E_{k-2}})
\cong\res{D(p-1)}{E_{k-2}}.
\]
The matrix in \cite{kjw}, after setting the last two coordinates equal to zero, gives
\[
\rank S'=(k-2)(p-1),\qquad \rank(S'^{p-1})=k-2.
\]
Hence the maximal Jordan type of $\res{D(1)}{E_{k-2}}$ is $[p]^{k-2}[1]^{2p-2}$. By Proposition \ref{p:wedgeJT}, the maximal Jordan type of $\res{D(p-1)}{E_{k-2}}$ is its $(p-1)$st exterior power. Since $\bigwedge^{p-1}([1]^{2p-2})$ contains a trivial summand, this exterior power contains $[1]$. It is therefore not free, and neither is $\res{D(p-1)}{\langle u_\alpha\rangle}$.
\end{proof}

\begin{thm}\label{t.main}
$\rk{E_k}{D(p-1)}=V(p_k)$.
\end{thm}

\begin{proof}
By Corollary~\ref{lem:complement} and Theorem~\ref{t.umaxDr},
\[
\rk{E_k}{D(p-1)}=\U{E_k}{D(p-1)}^c
\supseteq V(f_k)^c\cap V(p_k).
\]
On the other hand, Lemma \ref{l.inclusion} gives
\[\rk{E_k}{D(p-1)}\supseteq \bigcup_{1\leq i\neq j\leq k}V(x_i,x_j)=V(f_k)\cap V(p_k).\]
These two inclusions show that $\rk{E_k}{D(p-1)}\supseteq V(p_k)$. The reverse inclusion is \cite[Lemma 5.2]{kjw}.
\end{proof}

The arguments of \cite[Corollaries 5.6 and 5.7]{kjw} now give the following consequences.
\begin{cor}
$\rk{E_k}{D(kp-p-1)}=V(p_k)$.
\end{cor}
\begin{cor}
The complexities of the simple modules $D(p-1)$ and $D(kp-p-1)$ are $k-1$.
\end{cor}
\noindent
\textbf{Acknowledgments.}\par
The author acknowledges support from EPSRC grant EP/X035328/1.

\section*{Appendix A. Proof of Lemma \ref{l.lrprop2}}
\begingroup
\tikzset{every picture/.append style={scale=.90,transform shape,line cap=round,line join=round}}
\begin{proof}[Proof of Lemma \ref{l.lrprop2}]
We first prove part (i). For integers $m$ and $x\leq y$, we use the following notation:
\[
 \begin{tikzpicture}
 \matrix(m)[matrix of math nodes,nodes in empty cells]
 {&x&\\
 &&\\
  &\raisebox{0pt}[0.6\height][0\height]{$ \vdots $}&\\
   &&\\
  & y&\\};
\draw (m-5-1.south west) rectangle ([shift={(0,0.2)}]m-1-3.north east);
\draw[thick,decorate,decoration = {brace,mirror,raise=3pt}] (m-5-1.south west) -- (m-5-3.south east)
  node[pos=0.5, ,above=-20pt] {$m$};
    \end{tikzpicture} 
    \begin{tikzpicture}
 \matrix(m)[matrix of math nodes,nodes in empty cells]
 {x& x &\dots &x\\ 
 x+1 &x+1 &\dots & x+1 \\
\raisebox{0pt}[0.6\height][0\height]{$ \vdots $} &\raisebox{0pt}[0.6\height][0\height]{$ \vdots $} &\raisebox{0pt}[0.6\height][0\height]{$ \ddots $} &\raisebox{0pt}[0.6\height][0\height]{$ \vdots $}\\ 
 y& y&\dots & y\\};
  \draw ([shift={(-0.3,0)}]m-4-1.south west) rectangle ([shift={(0.22,0)}]m-1-4.north east);
   \node at (-2.8,0) {$=$};
   \draw[thick,decorate,decoration = {brace,mirror,raise=3pt}] ([shift={(-0.25,0)}]m-4-1.south west) -- ([shift={(0.22,0)}]m-4-4.south east)
  node[pos=0.5, ,above=-20pt] {$m$};
  \draw[thick,decorate,decoration = {brace,raise=3pt}] ([shift={(0.22,0)}]m-1-4.north east) -- ([shift={(0.22,0)}]m-4-4.south east)
  node[pos=0.5, ,right=4pt] {$y-x+1$};
    \end{tikzpicture}.
\]
Suppose that $\mu=(p^b,1)$ and $\beta=p-\bar\xi$ for some $\xi\vdash m$ with $\xi\neq (1^m)$; thus $\ell(\xi)=l\leq m-1$. Write $\xi=(\xi_1^{n_1},\ldots,\xi_t^{n_t})$, where $\xi_1>\cdots>\xi_t>0$, $n_1,\ldots,n_t>0$, and $\sum_{i=1}^t n_i=l$. Then
\begin{align*}
\beta&=p-\bar\xi\\
              &=p-((\sum_{i\leq t}n_i)^{\xi_t},(\sum_{i\leq t-1}n_i)^{\xi_{t-1}-\xi_t},\ldots,n_1^{\xi_1-\xi_2})\\
              &=((p-n_1)^{\xi_1-\xi_2},\ldots, (p-\sum_{i\leq t-1}n_i)^{\xi_{t-1}-\xi_t},(p-l)^{\xi_t}).
\end{align*}
It follows that $\bar\beta=(\xi_1^{p-l},(\xi_1-\xi_t)^{n_t},\ldots,(\xi_1-\xi_2)^{n_2})$. Let $T=T(A)$ be the tableau of shape $\bar\mu=(b+1,b^{p-1})$ corresponding to $A$, with entries $0,1,\ldots,r$, where $r=\ell(\bar\beta)=p-n_1$. Recall that $T^h$ contains the entries at most $h$ and has shape $\overline{a^h}$; in particular, $T^r=T$ and $T^0$ has shape $\bar\lambda$. Since $T_{b+1,b+1}(1)\geq T_{b+1,b+1}(h)$ for $2\leq h\leq r$, the entry in node $(1,b+1)$ is $q\in\{0,1\}$:

\[
 \begin{tikzpicture}
  \matrix (m)[
    matrix of math nodes,
    nodes in empty cells,
  ] {&   &  &   &      &   &   &  &   &      &  &   &  &   &      {\scriptstyle q}    \\
      &   &  &   &      &   &   &  &   &      &  &   &  &   &           \\
      &   &  &   &      &   &   &  &   &      &  &   &  &   &            \\
     &   &  &   &      &   &   &  &   &      &  &   &  &   &            \\
    &   &  &   &      &   &   &  &   &      &  &   &  &   &            \\
   &   &  &   &      &   &   &  &   &      &  &   &  &   &            \\
    &   &  &   &      &   &   &  &   &      &  &   &  &   &            \\
    &   &  &   &      &   &   &  &   &      &  &   &  &   &            \\
   &   &  &   &      &   &   &  &   &      &  &   &  &   &            \\
  };

  \draw (m-9-1.south west) rectangle (m-1-15.north west);
  \draw (m-1-15.south west) rectangle(m-1-15.north east);
  \draw[color=red] (m-9-10.south west) rectangle(m-9-14.north east);
   \draw[color=blue] ([shift={(0,0.18)}]m-1-1.north west) --(m-9-1.south west)--(m-9-9.south east)--(m-8-9.south east)--(m-8-14.south east)--(m-1-15.south west)--(m-1-15.south east)--(m-1-15.north east)--([shift={(0,0.18)}]m-1-1.north west);
  \pic{leftbrace={1,9,p}};
  \pic{overbrace={1,14,b}};
\draw[red,thick,decorate,decoration = {brace,mirror,raise=3pt}] (m-9-10.south west) -- (m-9-14.south east)
  node[pos=0.5, ,above=-20pt,red] {$\xi_1-\xi_2$};  
  \node at (-3.5,0) {$T=$};
  \node at (m-4-8) {\Large$\color{blue}T^{r-1}$};
\end{tikzpicture},\]
Since $T^r\backslash T^{r-1}$ is a removable horizontal strip of size $\xi_1-\xi_2$, filled with $p-n_1$, it must occupy the red box; the remaining tableau $T^{r-1}$, shown in blue, has shape $\overline{a^{r-1}}$. For $i=p-n_1-n_2+1,\ldots,p-n_1-1$,
\[
\bar\beta_i=\xi_1-\xi_2=T_{1,b+1}(i)\geq T_{b+1-(\xi_1-\xi_2),b}(i).
\]
Conversely,
\begin{align*}
 T_{b+1-(\xi_1-\xi_2),b}(i) &=T_{b+1-(\xi_1-\xi_2),b+1}(i)\\
&\geq T_{b+1-(\xi_1-\xi_2),b+1}(p-n_1)\\
&=\xi_1-\xi_2.
\end{align*}
Hence $T_{b+1-(\xi_1-\xi_2),b}(i)=\xi_1-\xi_2$, and each $T^i\backslash T^{i-1}$ is a horizontal strip of length $\xi_1-\xi_2$ stacked above $T\backslash T^{r-1}$. Thus $T\backslash T^{p-n_1-n_2}$ is an $n_2\times(\xi_1-\xi_2)$ rectangle:
\[\begin{tikzpicture}
 \node at (-1.6,1) {$T\backslash T^{p-n_1-n_2}=$};
 \draw (0,0) -- ++(0,1.7) --++(1.8,0)--++(0,-1.7)--++(-1.8,0);
 \node at (0.9,0.2)  {$\scriptscriptstyle p-n_1$};
  \node at (0.9,1)  {$\scriptstyle \vdots$};
 \node at (0.9,1.5)  {$\scriptstyle p-n_1-n_2+1$};
 \draw[decoration={brace,mirror,raise=2.5pt},decorate]
  (0,0) -- node[below=3pt] {$\scriptstyle\xi_1-\xi_2$} (1.8,0);
 \end{tikzpicture} \]
and \[\begin{tikzpicture}
\node at (-1.2,1) {$T^{p-n_1-n_2}=$};
 \draw (0,0) -- ++(0,2.5) --++(4,0)--++(0,-0.3)--++(-0.3,0)--++(0,-1)--++(-1,0)--++(0,-1.2)--++(-2.7,0);
  \draw[red] (0.01,0.01) -- ++(0,0.3)--++(2.68,0)--++(0,-0.3)--++(-2.68,0);
  \draw[red] (2.71,1.22)--++(0,0.3)--++(0.98,0)--++(0,-0.3)--++(-0.98,0);
   \draw[decoration={brace,mirror,raise=2.5pt},decorate]
  (2.7,1.2) -- node[below=3pt] {$\scriptstyle\xi_1-\xi_2$} (3.7,1.2);
 \end{tikzpicture}.\]
If $p-n_1-n_2>1$, the horizontal strip $T^{p-n_1-n_2}\backslash T^{p-n_1-n_2-1}$ does not contain node $(1,b+1)$ and must lie in the red strip above. Since $\bar\beta_{p-n_1-n_2}=\xi_1-\xi_3=T_{1,b+1}(p-n_1-n_2)$ and
\[
T_{b+1-(\xi_1-\xi_2),b+1}(p-n_1-n_2)
=T_{b+1-(\xi_1-\xi_2),b}(p-n_1-n_2)\geq\xi_1-\xi_2,
\]
we obtain
\[
T_{b+1-(\xi_1-\xi_2),b}(p-n_1-n_2)=\xi_1-\xi_2,
\qquad
T_{1,b-(\xi_1-\xi_2)}(p-n_1-n_2)=\xi_2-\xi_3.
\]
Thus
\[\begin{tikzpicture}
 \node at (-4,1) {$T^{p-n_1-n_2}\backslash T^{p-n_1-n_2-1}= $};
 \draw (0,0) -- ++(0,1.7) --++(1.5,0)--++(0,-1.7)--++(-2.9,0)--++(0,0.3)--++(1.4,0);
 \node at (0.75,0.15)  {$\scriptstyle p-n_1$};
  \node at (0.75,1)  {$\scriptstyle \vdots$};
 \node at (0.75,1.5)  {$\scriptstyle p-n_1-n_2$};
  \node at (-0.65,0.15)  {$\scriptstyle p-n_1-n_2$};
   \draw[decoration={brace,mirror,raise=2.5pt},decorate]
  (0,0) -- node[below=3pt] {$\scriptstyle\xi_1-\xi_2$} (1.5,0);
   \draw[decoration={brace,mirror,raise=2.5pt},decorate]
  (-1.4,0) -- node[below=3pt] {$\scriptstyle\xi_2-\xi_3$} (0,0);

 \end{tikzpicture};\]
 \[ \begin{tikzpicture}
 \node at (-1.4,1.5) {$T^{p-n_1-n_2-1}=$};
 \draw (0,0) -- ++(0,2.5) --++(4,0)--++(0,-0.3)--++(-0.3,0)--++(0,-1)--++(-1,0)--++(0,-0.9)--++(-0.7,0)--++(0,-0.3)--++(-2,0);
  \draw[red] (0.01,0.01) -- ++(0,0.3)--++(1.98,0)--++(0,-0.3)--++(-1.98,0);
  \draw[red] (2.01,0.31)--++(0,0.3)--++(0.68,0)--++(0,-0.3)--++(-0.68,0);
    \draw[red] (2.71,1.21)--++(0,0.3)--++(0.98,0)--++(0,-0.3)--++(-0.98,0);
     \draw[decoration={brace,mirror,raise=2.5pt},decorate]
  (2.7,1.2) -- node[below=3pt] {$\scriptstyle\xi_1-\xi_2$} (3.7,1.2);
   \draw[decoration={brace,mirror,raise=2.5pt},decorate]
  (2,0.3) -- node[below=3pt] {$\scriptstyle\xi_2-\xi_3$} (2.7,0.3);
 \end{tikzpicture}.\]
Moreover, $T^{p-n_1-n_2-1}\backslash T^{p-n_1-n_2-2}$ must lie in the red horizontal strip above. For $i=p-n_1-n_2-n_3+1,\ldots,p-n_1-n_2-1$, we have $\bar\beta_i=\xi_1-\xi_3=T_{1,b+1}(i)$ and
\begin{align*}
T_{1,b+1}(i)&\geq T_{b+1-(\xi_1-\xi_3),b}(i)
=T_{b+1-(\xi_1-\xi_3),b+1}(i)\\
&\geq T_{b+1-(\xi_1-\xi_3),b}(p-n_1-n_2)\\
&=\xi_1-\xi_3.
\end{align*}
Thus $T_{b+1-(\xi_1-\xi_3),b}(i)=\xi_1-\xi_3$. Similarly,
\[
T_{b+1-(\xi_1-\xi_2),b}(i)=\xi_1-\xi_2,
\qquad
T_{b+1-(\xi_1-\xi_3),b-(\xi_1-\xi_2)}(i)=\xi_2-\xi_3.
\]
Hence, for $i=p-n_1-n_2-n_3+1,\ldots,p-n_1-n_2-1$, each $T^i\backslash T^{i-1}$ is a horizontal strip consisting of two segments stacked above $T^{p-n_1-n_2}\backslash T^{p-n_1-n_2-1}$:
\[ \begin{tikzpicture}
 \node at (-3,1) {$T\backslash T^{i-1}=$};
 \draw (0,0) -- ++(0,1.7) --++(1.5,0)--++(0,-1.7)--++(-2.9,0)--++(0,1)--++(1.4,0);
 \node at (0.75,0.15)  {$\scriptstyle p-n_1$};
  \node at (0.75,1)  {$\scriptstyle \vdots$};
 \node at (0.75,1.6)  {$\scriptstyle i$};
  \node at (-0.65,0.15)  {$\scriptstyle p-n_1-n_2$};
  \node at (-0.65,0.9)  {$\scriptstyle i$};
   \draw[decoration={brace,mirror,raise=2.5pt},decorate]
  (0,0) -- node[below=3pt] {$\scriptstyle\xi_1-\xi_2$} (1.5,0);
   \draw[decoration={brace,mirror,raise=2.5pt},decorate]
  (-1.4,0) -- node[below=3pt] {$\scriptscriptstyle\xi_2-\xi_3$} (0,0);
 \end{tikzpicture},\] 
 \[\begin{tikzpicture}
 \node at (-1,1.5) {$T^{i-1}= $};
 \draw (0,0) -- ++(0,2.5) --++(4,0)--++(0,-0.3)--++(-0.3,0)--++(0,-1)--++(-1,0)--++(0,-0.5)--++(-0.7,0)--++(0,-0.7)--++(-2,0);
  \draw[red] (0.01,0.01) -- ++(0,0.3)--++(1.98,0)--++(0,-0.3)--++(-1.98,0);
    \draw[decoration={brace,mirror,raise=2.5pt},decorate]
  (2.7,1.2) -- node[below=3pt] {$\scriptstyle\xi_1-\xi_2$} (3.7,1.2);
   \draw[decoration={brace,mirror,raise=2.5pt},decorate]
  (2,0.7) -- node[below=3pt] {$\scriptstyle\xi_2-\xi_3$} (2.7,0.7);
 \end{tikzpicture}.\]

Repeating this process gives
 \[
 \begin{tikzpicture}
 \draw (0,0)--++(1.2,0)--++(0,1.5)--++(-1.2,0)--++(0,-1.5);
  \draw (1.2,0)--++(1.2,0)--++(0,2.1)--++(-1.2,0)--++(0,-2.1);
   \draw (3.9,0)--++(1.2,0)--++(0,2.7)--++(-1.2,0)--++(0,-2.7);
   \draw (5.1,0)--++(1.2,0)--++(0,3.3)--++(-1.2,0)--++(0,-3.3);
    \draw[blue] (0,0)--++(-1.5,0)--++(0,4.2)--++(8.2,0)--++(0,-0.4)--++(-0.4,0)--++(0,-0.5);
    \node at(-3,2.1) {$T=$};
    \draw[decoration={brace,raise=3.5pt},decorate]
  (-1.5,4.2) -- node[above=4pt] {$\displaystyle b+1$} (6.7,4.2);
   \draw[decoration={brace,raise=3.5pt},decorate]
  (-1.5,0) -- node[left=4pt] {$\displaystyle p$} (-1.5,4.2);
  \node[blue] at(2.2,3) {$T^1$};
  
  \node at(0.6,1.3){$\scriptscriptstyle2$};
  \node at (0.6,0.2){\scalebox{0.5}{$p-\sum_{i=1}^{t}n_i$}};
   \draw[decoration={brace,mirror,raise=2.5pt},decorate]
  (0,0) -- node[below=3pt] {$\scriptscriptstyle\xi_t$} (1.2,0);
  
   \node at(1.8,1.9){$\scriptscriptstyle2$};
  \node at (1.8,0.2){\scalebox{0.5}{$p-\sum_{i=1}^{t-1}n_i$}};
   \draw[decoration={brace,mirror,raise=2.5pt},decorate]
  (1.2,0) -- node[below=3pt] {$\scriptscriptstyle\xi_{t-1}-\xi_t$} (2.4,0);
  
  \node at(3.15,1.6){$\cdots$};
  
 \node at(4.5,2.5){$\scriptscriptstyle2$};
  \node at (4.5,0.2){\scalebox{0.5}{$p-n_1-n_2$}};
   \draw[decoration={brace,mirror,raise=2.5pt},decorate]
  (3.9,0) -- node[below=3pt] {$\scriptscriptstyle\xi_2-\xi_3$} (5.1,0);
  
   \node at(5.7,3.1){$\scriptscriptstyle2$};
  \node at (5.7,0.2){\scalebox{0.5}{$p-n_1$}};
   \draw[decoration={brace,mirror,raise=2.5pt},decorate]
  (5.1,0) -- node[below=3pt] {$\scriptscriptstyle\xi_1-\xi_2$} (6.3,0);
  
     \draw[decoration={brace,mirror,raise=2.5pt},decorate]
  (6.3,3.8) -- node[below=3pt] {$\scriptscriptstyle1$} (6.7,3.8);
 \end{tikzpicture}
 \]

If $T(1,b+1)=0$, the preceding argument gives
 \[
 \begin{tikzpicture}
 \draw (0,0)--++(1.2,0)--++(0,1.5)--++(-1.2,0)--++(0,-1.5);
  \draw (1.2,0)--++(1.2,0)--++(0,2.1)--++(-1.2,0)--++(0,-2.1);
   \draw (3.9,0)--++(1.2,0)--++(0,2.7)--++(-1.2,0)--++(0,-2.7);
   \draw (5.1,0)--++(1.2,0)--++(0,3.3)--++(-1.2,0)--++(0,-3.3);
    \draw[blue] (0,0)--++(-1.5,0)--++(0,4.2)--++(8.2,0)--++(0,-0.4)--++(-0.4,0)--++(0,-0.5);
    \node at(-3,2.1) {$T=$};
    \draw[decoration={brace,raise=3.5pt},decorate]
  (-1.5,4.2) -- node[above=4pt] {$\displaystyle b+1$} (6.7,4.2);
   \draw[decoration={brace,raise=3.5pt},decorate]
  (-1.5,0) -- node[left=4pt] {$\displaystyle p$} (-1.5,4.2);
  
  \node[blue] at(2.2,3) {$T^0$};
  
  \node at(0.6,1.3){$\scriptscriptstyle1$};
  \node at (0.6,0.2){\scalebox{0.5}{$p-\sum_{i=1}^{t}n_i$}};
   \draw[decoration={brace,mirror,raise=2.5pt},decorate]
  (0,0) -- node[below=3pt] {$\scriptscriptstyle\xi_t$} (1.2,0);
  
   \node at(1.8,1.9){$\scriptscriptstyle1$};
  \node at (1.8,0.2){\scalebox{0.5}{$p-\sum_{i=1}^{t-1}n_i$}};
   \draw[decoration={brace,mirror,raise=2.5pt},decorate]
  (1.2,0) -- node[below=3pt] {$\scriptscriptstyle\xi_{t-1}-\xi_t$} (2.4,0);
  
  \node at(3.15,1.6){$\cdots$};
  
 \node at(4.5,2.5){$\scriptscriptstyle1$};
  \node at (4.5,0.2){\scalebox{0.5}{$p-n_1-n_2$}};
   \draw[decoration={brace,mirror,raise=2.5pt},decorate]
  (3.9,0) -- node[below=3pt] {$\scriptscriptstyle\xi_2-\xi_3$} (5.1,0);
  
   \node at(5.7,3.1){$\scriptscriptstyle1$};
  \node at (5.7,0.2){\scalebox{0.5}{$p-n_1$}};
   \draw[decoration={brace,mirror,raise=2.5pt},decorate]
  (5.1,0) -- node[below=3pt] {$\scriptscriptstyle\xi_1-\xi_2$} (6.3,0);
  
     \draw[decoration={brace,mirror,raise=2.5pt},decorate]
  (6.3,3.8) -- node[below=3pt] {$\scriptscriptstyle1$} (6.7,3.8);
 \end{tikzpicture}
 \]
Since $T^0$ has shape $\overline{a^0}$, we obtain $\lambda=a^0=(p^h)\sqcup\delta$, where
\[\delta = (\underbrace{\sum_{i=1}^tn_i,\ldots,\sum_{i=1}^tn_i}_{\xi_t}, \underbrace{\sum_{i=1}^{t-1}n_i,\ldots, \sum_{i=1}^{t-1}n_i}_{\xi_{t-1}-\xi_t},\ldots, \underbrace{ n_1,\ldots, n_1}_{\xi_1-\xi_2}, 1).
\]
Otherwise, $q=1$ and $T^1\backslash T^0$ contains node $(1,b+1)$. In this case,
\[
T_{b+1-\xi_1,b}(1)=T_{b+1-\xi_1,b+1}(1)-1=\xi_1-1.
\]
Thus there is a column $j$, with $b+1-\xi_1\leq j\leq b$, such that $T_{i,i}(1)=1$ for $b+1-\xi_1\leq i\leq b$, except that $T_{j,j}(1)=0$. Hence
\[
\begin{tikzpicture}
 \node at (-2,1.5) {$T\backslash T^0= $};
 \draw (0,0) -- ++(0,1.6) --++(1.5,0)--++(0,-1.6)--++(-1.5,0);
 \node at (0.75,0.2)  {$\scriptscriptstyle p-\sum_{i=1}^tn_i$};
  \node at (0.75,0.9)  {$\scriptscriptstyle \vdots$};
 \node at (0.75,1.4)  {$\scriptscriptstyle 1$};
 \draw[decoration={brace,mirror,raise=2.5pt},decorate]
  (0,0) -- node[below=3pt] {$\scriptscriptstyle\xi_t$} (1.5,0);
   \node at (2.25,1.2)  {$\cdots$};
   \draw (3,0) -- ++(0,2)--++(0.3,0)--++(0,0.3) --++(1.2,0)--++(0,-2.3)--++(-1.5,0);
 \node at (3.75,0.2)  {$\scriptscriptstyle p-\sum_{i=1}^kn_i$};
  \node at (3.75,1.25)  {$\scriptscriptstyle \vdots$};
 \node at (3.75,2.1)  {$\scriptscriptstyle 1$};
 \draw[decoration={brace,mirror,raise=2.5pt},decorate]
  (3,0) -- node[below=3pt] {$\scriptscriptstyle\xi_k-\xi_{k+1}$} (4.5,0);
   \node at (5.25,1.2)  {$\cdots$};
  \draw (6,0) -- ++(0,2.8) --++(1.5,0)--++(0,-2.8)--++(-1.5,0);
 \node at (6.75,0.2)  {$\scriptscriptstyle p-n_1$};
  \node at (6.75,1.5)  {$\scriptscriptstyle \vdots$};
 \node at (6.75,2.6)  {$\scriptscriptstyle 1$};
 \draw[decoration={brace,mirror,raise=2.5pt},decorate]
  (6,0) -- node[below=3pt] {$\scriptscriptstyle\xi_1-\xi_2$} (7.5,0);
  \draw (7.5,3.3)--++(0.4,0)--++(0,0.4)--++(-0.4,0)--++(0,-0.4);
  \node at (7.7,3.5) {$\scriptscriptstyle 1$};
 \end{tikzpicture}.\] 
In this case, we obtain that $\lambda=a^0=(p^h)\sqcup \delta$ where \[\delta = (\underbrace{\sum_{i=1}^tn_i,\ldots,\sum_{i=1}^tn_i}_{\xi_t}, \ldots,\underbrace{\sum_{i=1}^{k}n_i+1,\sum_{i=1}^{k}n_i,\ldots, \sum_{i=1}^{k}n_i}_{\xi_{k}-\xi_{k+1}},\ldots, \underbrace{ n_1,\ldots, n_1}_{\xi_1-\xi_2}, 1).
\]
In both cases, $\delta$ is a partition of $m+1$ with at least two parts, so $\delta\neq(m+1)$. This proves part (i).
For part (ii), suppose that $\mu=(p^b,p-1)$ and $(m)\neq\beta\vdash m$. Write $\bar\beta=(\xi_1^{n_1},\ldots,\xi_t^{n_t})$, where $\xi_1>\cdots>\xi_t>0$, $n_1,\ldots,n_t>0$, and $\sum_{i=1}^t n_i=r=\ell(\bar\beta)<m$. Let $T=T(A)$ be the corresponding tableau of shape $\bar\mu=((b+1)^{p-1},b)$, with entries $0,1,\ldots,r$:
\[ \begin{tikzpicture}
\draw (0,0) -- ++(0,2.5) --++(4,0)--++(0,-2.1)--++(-0.4,0)--++(0,-0.4)--++(-3.6,0);
  \draw[red] (0.01,0.01) -- ++(0,0.38)--++(3.58,0)--++(0,-0.38)--++(-3.58,0);
  \draw[red] (3.61,0.41)--++(0,0.38)--++(0.38,0)--++(0,-0.38)--++(-0.38,0);
   \draw[decoration={brace,raise=3pt},decorate]
  (0,2.5) -- node[above=3.5pt] {$\displaystyle b+1$} (4,2.5);
   \draw[decoration={brace,raise=3pt},decorate]
  (0,0) -- node[left=3.5pt] {$\displaystyle p$} (0,2.5);
  \node at (-1,1.2) {$T=$};
 \end{tikzpicture}.\]
Since $|T\backslash T^{r-1}|=\bar\beta_r=\xi_t$, the difference $T\backslash T^{r-1}$ is a removable horizontal strip of size $\xi_t$ contained in the red strip. If it contains node $(p-1,b+1)$, then
\[
\begin{tikzpicture}
\node at (-2,.5) {$T\backslash T^{r-1}=$};
\draw (0,0) -- ++(0,0.4) --++(1.6,0)--++(0,0.4)--++(-0.4,0)--++(0,-0.8)--++(-1.2,0);
 \draw[decoration={brace,mirror,raise=2.5pt},decorate]
  (0,0) -- node[below=3pt] {$\scriptscriptstyle\xi_t-1$} (1.2,0);
  \draw[decoration={brace,mirror,raise=2.5pt},decorate]
  (1.2,0.4) -- node[below=3pt] {$\scriptscriptstyle1$} (1.6,0.4);
  \node at(0.6,0.2) {$\scriptstyle r$};
  \node at(1.4,0.6) {$\scriptstyle r$};
\end{tikzpicture}\]
\[\begin{tikzpicture}
\node at (-1,1) {$T^{r-1}= $};
\draw (0,0) -- ++(0,1.8) --++(2.5,0)--++(0,-1.2)--++(-0.3,0)--++(0,-0.3)--++(-0.8,0)--++(0,-0.3)--++(-1.4,0);
   \draw[decoration={brace,mirror,raise=2pt},decorate]
  (1.4,0.3) -- node[below=2.5pt] {$\scriptscriptstyle\xi_t-1$} (2.2,0.3);
  \draw[decoration={brace,mirror,raise=2pt},decorate]
  (2.2,0.6) -- node[below=2.5pt] {$\scriptscriptstyle 1$} (2.5,0.6);
 \end{tikzpicture}.\]
The argument for part (i) now gives
 \[\begin{tikzpicture}
 \node at (-1.5,1.5) {$T\backslash T^0=$};
 \draw (0,0)--++(0,1.3)--++(1.2,0)--++(0,-1.3)--++(-1.2,0);
   \node at(0.6,0.2) {$\scriptscriptstyle n_1$};
  \node at(0.6,1.1) {$\scriptscriptstyle 1$};
  \node at (0.6, 0.75){$\vdots$};
    \node at(1.6,1) {$\cdots$};
 \draw (2,0)--++(0,1.7)--++(1.2,0)--++(0,-1.7)--++(-1.2,0);
 \node at(2.6,0.2) {$\scriptscriptstyle \sum_{i=1}^{t-1}n_i$};
  \node at(2.6,1.5) {$\scriptscriptstyle 1$};
  \node at(2.6,0.95) {$\vdots$};
  
  \draw (3.2,0)--++(0,2.3)--++(1,0)--++(0,-2.3)--++(-1,0);
      \node at(3.7,0.2) {$\scriptscriptstyle r$};
   \node at(3.7,2.1) {$\scriptscriptstyle 1$};
   \node at(3.7,1.2) {$\vdots$};
   
  \draw (4.2,0.4)--++(0,2.3)--++(0.4,0)--++(0,-2.3)--++(-0.4,0);  
   \node at(4.4,0.6) {$\scriptscriptstyle r$};
   \node at(4.4,2.5) {$\scriptscriptstyle 1$};
   \node at(4.4,1.6) {$\vdots$};
   \draw[decoration={brace,mirror,raise=2pt},decorate]
  (0,0) -- node[below=2.5pt] {$\scriptscriptstyle\xi_1-\xi_2$} (1.2,0);
  \draw[decoration={brace,mirror,raise=2pt},decorate]
  (2,0) -- node[below=2.5pt] {$\scriptscriptstyle \xi_{t-1}-\xi_t$} (3.2,0); 
  \draw[decoration={brace,mirror,raise=2pt},decorate]
  (3.2,0) -- node[below=2.5pt] {$\scriptscriptstyle \xi_{t}-1$} (4.2,0); 
  \draw[decoration={brace,mirror,raise=2pt},decorate]
  (4.2,0.4) -- node[below=2.5pt] {$\scriptscriptstyle 1$} (4.6,0.4);
  \draw[decoration={brace,raise=2pt},decorate]
  (4.2,2.3) -- node[left=2.5pt] {$\scriptscriptstyle 1$} (4.2,2.7);    
 \end{tikzpicture},\]
which lies in the lower-right corner of $T$. Reading this diagram gives $\lambda=a^0=(p^h)\sqcup(p-\delta)$, where $h=b+1-\xi_1$ and
\[\delta=(r+1, \underbrace{r,\ldots, r}_{\xi_t-1},\underbrace{ \sum_{i=1}^{t-1}n_i,\ldots,\sum_{i=1}^{t-1}n_i}_{\xi_{t-1}-\xi_t},\ldots, \underbrace{n_1,\ldots,n_1}_{\xi_1-\xi_2}).
\]
Here $|\delta|=m+1$ and, since $r<m$, $\delta_1<m+1$; hence $\delta\neq(m+1)$. It remains to consider the case in which $T\backslash T^{r-1}$ does not contain node $(p-1,b+1)$. Then
\[\begin{tikzpicture}
\node at (-2,0.2) {$T\backslash T^{r-1}= $};
\draw (0,0) -- ++(0,0.4) --++(1.2,0)--++(0,-0.4)--++(-1.2,0);
 \draw[decoration={brace,mirror,raise=2.5pt},decorate]
  (0,0) -- node[below=3pt] {$\scriptscriptstyle\xi_t$} (1.2,0);
   \node at(0.6,0.2) {$\scriptstyle r$};
\end{tikzpicture},\]
and
\[\begin{tikzpicture}
\node at (-1,1) {$T^{r-1}=$};
\draw (0,0) -- ++(0,1.8) --++(2.5,0)--++(0,-1.5)--++(-0.9,0)--++(0,-0.3)--++(-1.6,0);
   \draw[decoration={brace,mirror,raise=2pt},decorate]
  (1.6,0.3) -- node[below=2.5pt] {$\scriptscriptstyle\xi_t+1$} (2.5,0.3);
  \draw[red] (1.6,0.32)--++(0,0.3)--++(0.9,0)--++(0,-0.3)--++(-0.9,0); 
 \end{tikzpicture}. \]
Since $|T^{r-1}\backslash T^{r-2}|=\xi_t$ and $T_{b+1-\xi_t,b+1}(r-1)\geq T_{b+1-\xi_t,b+1}(r)=\xi_t$, the difference $T^{r-1}\backslash T^{r-2}$ is a horizontal strip of size $\xi_t$ contained in the red strip. Thus
\[ \begin{tikzpicture}
\node at (-2,0.5) {$T\backslash T^{r-2}= $};
\draw (0,0) -- ++(0,0.4) --++(1.6,0)--++(0,0.4)--++(-1.2,0)--++(0,-0.4);
\draw (0,0)--++(1.2,0)--++(0,0.4);
 \draw[decoration={brace,mirror,raise=2.5pt},decorate]
  (0,0) -- node[below=3pt] {$\scriptscriptstyle\xi_t$} (1.2,0);
 \draw[decoration={brace,raise=2.5pt},decorate]
  (0.4,0.8) -- node[above=3pt] {$\scriptscriptstyle\xi_t$} (1.6,0.8);
  \node at(0.6,0.2) {$\scriptstyle r$};
  \node at(1,0.6) {$\scriptstyle r-1$};
\end{tikzpicture}.\] 
Iterating this argument gives
 \[\begin{tikzpicture}
\node at (-2,1) {$T\backslash T^{r-n_t}=$};
\draw (0,0) -- ++(0,0.4) --++(1.2,0)--++(0,-0.4)--++(-1.2,0);
\draw (0.4,0.4)--++(1.2,0)--++(0,1.9)--++(-1.2,0)--++(0,-1.9);
 \draw[decoration={brace,mirror,raise=2.5pt},decorate]
  (0,0) -- node[below=3pt] {$\scriptscriptstyle\xi_t$} (1.2,0);
 \draw[decoration={brace,raise=2.5pt},decorate]
  (0.4,2.3) -- node[above=3pt] {$\scriptscriptstyle\xi_t$} (1.6,2.3);
  \node at(0.6,0.2) {$\scriptstyle r$};
  \node at(1,0.6) {$\scriptstyle r-1$};
  \node at(1,1.35){$\vdots$};
  \node at(1,2.1) {$\scriptstyle r-n_t+1$};
\end{tikzpicture}\] 
\[\begin{tikzpicture}
\node at (-1.5,1) {$T^{r-n_t}=$};
\draw (0,0) -- ++(0,2.2) --++(2.7,0)--++(0,-1.2)--++(-0.9,0)--++(0,-0.7)--++(-0.3,0)--++(0,-0.3)--++(-1.5,0);
   \draw[decoration={brace,mirror,raise=2pt},decorate]
  (1.8,1) -- node[below=2.5pt] {$\scriptscriptstyle\xi_t-1$} (2.7,1);
\draw[decoration={brace,mirror,raise=2pt},decorate]
  (1.5,0.3) -- node[below=2.5pt] {$\scriptscriptstyle 1$} (1.8,0.3);
 \end{tikzpicture}.\]
We now prove part (ii), for any $m$ and sufficiently large $b$, by induction on $t\geq2$. If $t=2$, then $\bar\beta=(\xi_1^{n_1},\xi_2^{n_2})$, and there are two possibilities:
\[\begin{tikzpicture}
\node at (-3,1) {$T\backslash T^{n_1-1}= $};
\draw (0,0) -- ++(0,0.4) --++(1.3,0)--++(0,-0.4)--++(-1.3,0);
\draw (0.4,0.4)--++(1.3,0)--++(0,1.9)--++(-1.3,0)--++(0,-1.9);
\draw(0,0)--++(-1.2,0)--++(0,0.4)--++(1.2,0);
 \draw[decoration={brace,mirror,raise=2.5pt},decorate]
  (0,0) -- node[below=3pt] {$\scriptscriptstyle\xi_2$} (1.3,0);
 \draw[decoration={brace,raise=2.5pt},decorate]
  (0.4,2.3) -- node[above=3pt] {$\scriptscriptstyle\xi_2$} (1.7,2.3);
   \draw[decoration={brace,mirror,raise=2.5pt},decorate]
  (-1.2,0) -- node[below=3pt] {$\scriptscriptstyle\xi_{1}-\xi_2$} (0,0);
 \draw[decoration={brace,raise=2.5pt},decorate]
  (0,0.4) -- node[above=3pt] {$\scriptscriptstyle1$} (0.4,0.4);
  \node at(0.65,0.2) {$\scriptscriptstyle n_1+n_2$};
  \node at(1.05,0.6) {$\scriptscriptstyle n_1+n_2-1$};
  \node at(1.05,1.35){$\vdots$};
  \node at(1.05,2.1) {$\scriptscriptstyle n_1$};
  \node at(-0.6,0.2) {$\scriptscriptstyle n_1$};
\end{tikzpicture}, \ \ \begin{tikzpicture}
\draw (0,0) -- ++(0,0.4) --++(1.3,0)--++(0,-0.4)--++(-1.3,0);
\draw (0.4,0.4)--++(1.3,0)--++(0,1.9)--++(-1.3,0)--++(0,-1.9);
\draw(0,0)--++(-1,0)--++(0,0.4)--++(1,0);
\draw (0,0.4)--++(0,0.4)--++(0.4,0);
 \draw[decoration={brace,mirror,raise=2.5pt},decorate]
  (0,0) -- node[below=3pt] {$\scriptscriptstyle\xi_2$} (1.3,0);
 \draw[decoration={brace,raise=2.5pt},decorate]
  (0.4,2.3) -- node[above=3pt] {$\scriptscriptstyle\xi_2$} (1.7,2.3);
   \draw[decoration={brace,mirror,raise=2.5pt},decorate]
  (-1,0) -- node[below=3pt] {$\scriptscriptstyle\xi_{1}-\xi_2-1$} (0,0);
 \draw[decoration={brace,raise=2.5pt},decorate]
  (0,0.8) -- node[above=3pt] {$\scriptscriptstyle1$} (0.4,0.8);
  \node at(0.65,0.2) {$\scriptscriptstyle n_1+n_2$};
  \node at(1.05,0.6) {$\scriptscriptstyle n_1+n_2-1$};
  \node at(1.05,1.35){$\vdots$};
  \node at(1.05,2.1) {$\scriptscriptstyle n_1$};
  \node at(-0.5,0.2) {$\scriptscriptstyle n_1$};
   \node at(0.2,0.6) {$\scriptscriptstyle n_1$};
\end{tikzpicture}.\]
Accordingly, $T\backslash T^0$ is one of the following:
\[
\begin{tikzpicture}
\draw (0,0) -- ++(0,0.4) --++(1.3,0)--++(0,-0.4)--++(-1.3,0);
\draw (0.4,0.4)--++(1.3,0)--++(0,1.9)--++(-1.3,0)--++(0,-1.9);
\draw(0,0)--++(-1.2,0)--++(0,0.4)--++(1.2,0);
\draw(-0.8,0.4)--++(0,1)--++(1.2,0);
 \draw[decoration={brace,mirror,raise=2.5pt},decorate]
  (0,0) -- node[below=3pt] {$\scriptscriptstyle\xi_2$} (1.3,0);
 \draw[decoration={brace,raise=2.5pt},decorate]
  (0.4,2.3) -- node[above=3pt] {$\scriptscriptstyle\xi_2$} (1.7,2.3);
   \draw[decoration={brace,mirror,raise=2.5pt},decorate]
  (-1.2,0) -- node[below=3pt] {$\scriptscriptstyle\xi_{1}-\xi_2$} (0,0);
   \draw[decoration={brace,raise=2.5pt},decorate]
  (-0.8,1.4) -- node[above=3pt] {$\scriptscriptstyle\xi_{1}-\xi_2$} (0.4,1.4);
  \draw[decoration={brace,mirror,raise=2.5pt},decorate]
  (1.3,0.4) -- node[below=3pt] {$\scriptscriptstyle 1$} (1.7,0.4);
  \node at(0.65,0.2) {$\scriptscriptstyle n_1+n_2$};
  \node at(1.05,0.6) {$\scriptscriptstyle n_1+n_2-1$};
  \node at(1.05,1.35){$\vdots$};
  \node at(1.05,2.1) {$\scriptscriptstyle n_1$};
  \node at(-0.6,0.2) {$\scriptscriptstyle n_1$};
   \node at(-0.2,0.6) {$\scriptscriptstyle n_1-1$};
  \node at(-0.2,1.2) {$\scriptscriptstyle 1$};
  \node at(-0.2,0.9) {$\scriptscriptstyle \vdots$};
\end{tikzpicture}\; \;\; \ \ 
\begin{tikzpicture}
\draw (0,0) -- ++(0,0.4) --++(1.3,0)--++(0,-0.4)--++(-1.3,0);
\draw (0.4,0.4)--++(1.3,0)--++(0,1.9)--++(-1.3,0)--++(0,-1.9);
\draw(0,0)--++(-1.2,0)--++(0,0.4);
\draw (0,0.4)--++(0,1.4)--++(0.4,0);
\draw(-1.2,0.4)--++(0,1)--++(1.2,0);
 \draw[decoration={brace,mirror,raise=2.5pt},decorate]
  (0,0) -- node[below=3pt] {$\scriptscriptstyle\xi_2$} (1.3,0);
 \draw[decoration={brace,raise=2.5pt},decorate]
  (0.4,2.3) -- node[above=3pt] {$\scriptscriptstyle\xi_2$} (1.7,2.3);
   \draw[decoration={brace,mirror,raise=2.5pt},decorate]
  (-1.2,0) -- node[below=3pt] {$\scriptscriptstyle\xi_{1}-\xi_2-1$} (0,0);
    \draw[decoration={brace,mirror,raise=2.5pt},decorate]
  (1.3,0.4) -- node[below=3pt] {$\scriptscriptstyle 1$} (1.7,0.4);
  \node at(0.65,0.2) {$\scriptscriptstyle n_1+n_2$};
  \node at(1.05,0.6) {$\scriptscriptstyle n_1+n_2-1$};
  \node at(1.05,1.35){$\vdots$};
  \node at(1.05,2.1) {$\scriptscriptstyle n_1$};
  \node at(-0.6,0.2) {$\scriptscriptstyle n_1$};
  \node at(-0.6,1.2) {$\scriptscriptstyle 1$};
  \node at(-0.6,0.8) {$\scriptscriptstyle \vdots$};
  \node at(0.2,0.6) {$\scriptscriptstyle n_1$};
  \node at(0.2,1.6) {$\scriptscriptstyle 1$};
  \node at(0.2,1.2){$\scriptscriptstyle \vdots$};
   \draw[decoration={brace,raise=2.5pt},decorate]
  (0,1.8) -- node[above=3pt] {$\scriptscriptstyle1$} (0.4,1.8);
\end{tikzpicture}.\]
Hence $\lambda=a^0=(p^h)\sqcup(p-\delta)$, where either
\[
\delta=((n_1+n_2)^{\xi_2},n_1^{\xi_1-\xi_2},1),\quad h=b-\xi_1,
\]
or
\[
\delta=((n_1+n_2)^{\xi_2},n_1+1,n_1^{\xi_1-\xi_2-1}),\quad h=b+1-\xi_1.
\]
Thus part (ii) holds for $t=2$. Now assume it holds for $2\leq t'<t$. For $1\leq i\leq r-1$,
\[
\xi_t\geq T_{b+2-\xi_t,b+1}(i)\geq T_{b+2-\xi_t,b+1}(r-1)=\xi_t,
\]
so $T_{b+2-\xi_t,b+1}(i)=\xi_t$. Hence the last $\xi_t$ columns of $T$ have the form
\[ \begin{tikzpicture}
\draw (0,0) -- ++(0,2.5) --++(4,0)--++(0,-2.1)--++(-0.4,0)--++(0,-0.4)--++(-3.6,0);
   \draw[decoration={brace,raise=3pt},decorate]
  (0,2.5) -- node[above=3.5pt] {$\displaystyle b+1$} (4,2.5);
  \draw (2.4,0)--++(0,0.4)--++(1.2,0);
  \node at (3,0.2){$\scriptscriptstyle r$};
  \draw(2.8,0.4)--++(0,1.3)--++(1.2,0);
  \node at (3.4,0.6){$\scriptscriptstyle r-1$};
  \node at (3.4,1.5){$\scriptscriptstyle 1$};
  \node at (3.4,1.05){$\scriptscriptstyle \vdots$};
  \draw[decoration={brace,mirror,raise=2.5pt},decorate]
  (2.4,0) -- node[below=3pt] {$\scriptscriptstyle\xi_{t}-1$} (3.6,0);
  \draw[decoration={brace,mirror,raise=2.5pt},decorate]
  (3.6,0.4) -- node[below=3pt] {$\scriptscriptstyle1$} (4,0.4);
  \draw[red] (2.38,-0.2)--++(0,0.62)--++(0.4,0)--++(0,2.32);
   \draw[decoration={brace,raise=3pt},decorate]
  (0,0) -- node[left=3.5pt] {$\displaystyle p$} (0,2.5);
  \node at (-1,1.2) {$T=$};
 \node at (1.2,1.2) {$S$};
 \end{tikzpicture}.\]
Cutting $T$ along the red line gives the tableau $S$ on the left. It has shape $\overline{(p^{b'},p-1)}$, where $b'=b-\xi_t$, and entries $1,\ldots,r'$, where $r'=\sum_{i=1}^{t-1}n_i$. Define $\lambda'$ by $\lambda=\lambda'\sqcup((p-r)^{\xi_t})$, let $\beta'$ be conjugate to $((\xi_1-\xi_t)^{n_1},\ldots,(\xi_{t-1}-\xi_t)^{n_{t-1}})$, and put $\mu'=(p^{b'},p-1)$. Let $A'=[\alpha^0,\ldots,\alpha^{r'}]$, where $a^i=\alpha^i\sqcup((p-r+i)^{\xi_t})$ for $0\leq i\leq r'$. Then $A'$ is an LR sequence of type $[\lambda',\beta';\mu']$, and $S$ is its corresponding tableau. Moreover, $|\beta'|=m-r\xi_t$ and $\ell(\beta')>1$ because $t-1\geq2$. By induction,
\[
\lambda'=(p^h)\sqcup(p-\delta'),\qquad
|\delta'|=m-r\xi_t+1,\qquad
\ell(\delta')>1,
\]
where $h=b'-\ell(\delta')+1$. Therefore
\[
\lambda=(p^h)\sqcup(p-\delta),\qquad
\delta=(r^{\xi_t})\sqcup\delta',
\]
with $|\delta|=m+1$ and $\ell(\delta)>1$. Finally,
\[
h=b'-\ell(\delta')+1=(b-\xi_t)-(\ell(\delta)-\xi_t)+1
=b-\ell(\delta)+1.
\]
This completes the proof.
\end{proof}
\endgroup

\end{document}